\documentclass[11pt,a4paper]{article}

\usepackage[T1]{fontenc}
\usepackage[utf8]{inputenc}
\usepackage{lmodern}
\usepackage{microtype}
\usepackage{amsmath,amssymb,amsthm,mathtools}
\usepackage{geometry}
\usepackage{enumitem}
\usepackage{xcolor}
\usepackage{hyperref}
\usepackage[nameinlink,noabbrev]{cleveref}

\hypersetup{colorlinks=true,linkcolor=blue,citecolor=blue,urlcolor=blue}
\numberwithin{equation}{section}

\newtheorem{theorem}{Theorem}[section]
\newtheorem{proposition}[theorem]{Proposition}
\newtheorem{lemma}[theorem]{Lemma}
\newtheorem{corollary}[theorem]{Corollary}
\theoremstyle{definition}

\newtheorem{example}[theorem]{Example}
\theoremstyle{remark}
\newtheorem{remark}[theorem]{Remark}

\newcommand{\R}{\mathbb R}
\newcommand{\Vol}{\operatorname{Vol}}
\DeclareMathOperator{\RelCap}{Cap}
\newcommand{\Per}{\operatorname{Per}}
\newcommand{\Ric}{\operatorname{Ric}}

\newcommand{\dd}{\,\mathrm d}
\newcommand{\1}{\mathbf 1}

\title{\textbf{An Isoperimetric Comparison for a Variational Threshold Associated with the Riemannian \(p\)-Laplacian}}
\author{
\begin{tabular}{c}
Ababacar Sadikhe DJIT\'E\\
\small \texttt{ababacarsadikhe.djite@ucad.edu.sn}\\[1.5mm]
Diaraf SECK\\
\small \texttt{diaraf.seck@ucad.edu.sn}\\[2mm]
\small Ecole Doctorale de Math\'ematiques et Informatique\\
\small Universit\'e Cheikh Anta Diop de Dakar, S\'en\'egal\\
\small Laboratoire des Math\'ematiques de la D\'ecision et d'Analyse Num\'erique\\
\small BP 16889, Dakar Fann, S\'en\'egal
\end{tabular}
}
\date{}

\begin{document}
\maketitle

\begin{abstract}
%We study a variational threshold associated with a one-phase functional involving the Riemannian \(p\)-Dirichlet energy and a volume penalization. We first characterize the threshold as an infimum of a capacity-to-volume quotient. We then prove an isoperimetric comparison theorem: domination of the isoperimetric profile of a domain by that of a model space implies a lower bound for its variational threshold in terms of the threshold on a model geodesic ball. For geodesic balls in simply connected space forms, we derive an explicit one-dimensional formula for the threshold and prove its strict monotonicity with respect to sectional curvature, both at fixed radius and at fixed volume. We also explain how the normalized L\'evy--Gromov inequality yields a comparison result on compact manifolds with a positive Ricci lower bound. Finally, we give a product-cylinder example showing why a nonnegative Ricci lower bound alone cannot be used to infer a global Euclidean isoperimetric comparison on a noncompact manifold.

We study a variational threshold $\Lambda_{p,M}(K)$ associated with a
one-phase free boundary functional involving the Riemannian $p$-Dirichlet energy and a volume
penalization on a bounded domain $K$ of a Riemannian manifold $M$. In the Euclidean case, this
threshold is the smallest parameter for which the interior Bernoulli problem has a solution
minimizing the energy; it bounds the Bernoulli constant from above, and the inequality is strict
for balls. We characterize $\Lambda_{p,M}(K)$ as the infimum of a capacity-to-volume quotient and
prove that this infimum is always attained. We then prove an isoperimetric comparison theorem:
domination of the isoperimetric profile of $M$ by that of a model space implies a lower bound for
$\Lambda_{p,M}(K)$ in terms of the threshold of a model geodesic ball. This extends a theorem of
Daners and Kawohl. For geodesic balls in space forms, we give an explicit one-dimensional formula
for the threshold and prove its strict monotonicity in the curvature, at fixed radius and at fixed
volume. As applications, we obtain comparison results on compact manifolds with positive Ricci
curvature, through the normalized Lévy--Gromov inequality, and on manifolds with nonnegative Ricci
curvature and Euclidean volume growth, where the factor $\mbox{AVR}^{1/N}$ appears; a product-cylinder
example shows that nonnegative Ricci curvature alone is not sufficient. Finally, we characterize
the equality cases: geodesic balls in space forms, round spheres and geodesic balls in the
Lévy--Gromov setting, Euclidean space and balls in the setting of Euclidean volume growth.
\end{abstract}

%\noindent\textbf{Keywords:} Riemannian \(p\)-Laplacian; variational threshold; relative \(p\)-capacity; isoperimetric profile; P\'olya--Szeg\H{o} inequality; L\'evy--Gromov inequality; space forms.

%\noindent\textbf{Mathematics Subject Classification (2020):} 49Q10, 35J92, 53C21.


%===========================================================
\section{Introduction}
%===========================================================

Variational problems involving free boundaries arise naturally in geometric analysis, nonlinear potential theory, and shape optimization. A basic class of such problems is described by one-phase functionals in which a Dirichlet energy competes with a volume penalization. The positivity region of a minimizer is not prescribed in advance, and its geometry is part of the variational problem.

Let \((M,g)\) be an \(N\)-dimensional Riemannian manifold and let \(K\subset M\) be a bounded domain. For \(1<p<\infty\) and \(\lambda>0\), we consider
\begin{equation}
\mathcal J_{\lambda,p}^{K}(u)
=
\int_K |\nabla_g u|_g^p\,\dd V_g
+(p-1)\lambda^p\Vol_g(\{u<1\}),
\qquad u\in W_0^{1,p}(K).
\label{eq:1-1}\end{equation}
The zero function is an admissible competitor. We define the variational threshold as the infimum of the parameters for which a nonzero global minimizer exists.

\paragraph{Relation with the Bernoulli free boundary problem.} Setting $v=1-u$, the functional
\eqref{eq:1-1} becomes
\[
\int_K|\nabla_g v|_g^p\dd V_g+(p-1)\lambda^p\Vol_g(\{v>0\}),\qquad v=1\ \text{on }\partial K,
\]
which is the one-phase functional of Alt and Caffarelli \cite{AC81} for $p=2$, and its
$p$-Laplacian version studied by Danielli and Petrosyan \cite{DP05}; see also \cite{BayramiFotouhi2023} for the
two-phase problem. In the Euclidean setting, nonzero global minimizers are therefore weak
solutions, in the sense of \cite{AC81,DP05}, of the interior Bernoulli problem
\[
\Delta_p u=0\ \text{in }K\setminus\overline\Omega,\qquad u=0\ \text{on }\partial K,\qquad u=1\ \text{and}\
|\nabla u|=\lambda\ \text{on }\partial\Omega,\qquad \Omega=\{u=1\}^\circ ,
\]
the normalization $(p-1)\lambda^p$ in \eqref{eq:1-1} being chosen so that the free boundary condition reads
$|\nabla u|=\lambda$. 
%The variational threshold $\Lambda_{p,M}(K)$ must not be confused with the
%\emph{Bernoulli constant} $c_K$, the infimum of the parameters $\lambda$ for which the interior
%problem admits a solution \cite{FR97,LySeck2004}. When $c_K$ is defined through weak solutions in the
%above sense, one has  $c_K\le\Lambda_{p,0}(K)$ (still noted $\Lambda_{p,\R^N}(K)$) %$c_K\le\Lambda_{p,\R^N}(K)$
%, and the inequality is strict for balls
%(Proposition~\ref{prop:ball} below).

For a domain $K\subset\R^N$ we write $\Lambda_{p,\R^N}(K)$ for $\Lambda_{p,M}(K)$ with $M=\R^N$;
on model balls we use the notation $\Lambda_{p,\kappa}(B_R^\kappa)$ of Section~4, so that
$\Lambda_{p,0}(B_R)=\Lambda_{p,\R^N}(B_R)$. The variational threshold must not be confused with
the Bernoulli constant $c_K$, the infimum of the parameters $\lambda$ for which the interior
problem admits a solution \cite{FR97,LySeck2004}. In the Euclidean setting, Proposition~7.1 shows that
for $\lambda=\Lambda_{p,\R^N}(K)$ the functional \eqref{eq:1-1} has a nonzero global minimizer, which is a
weak solution of the interior problem by \cite{AC81,DP05}. Hence $c_K\le\Lambda_{p,\R^N}(K)$ when
$c_K$ is defined through weak solutions, and the inequality is strict for balls
(Proposition~2.1).

Flucher and Rumpf \cite{FR97} asked whether $c_K\ge c_{K^*}$, where $K^*$ is the ball with the same
volume as $K$; see also \cite{LySeck2004}. The variational counterpart of this question was answered by
Daners and Kawohl \cite{DanersKawohl2010}: in $\R^N$, the minimal parameter for which the interior problem has
a solution minimizing the energy is minimal for the ball, and essentially only for the ball. 
Up to
the normalization of the parameter, this is the Euclidean case $M=\R^N$, $\kappa=0$, of
Theorem~3.4. The present paper extends this comparison to Riemannian manifolds under an
isoperimetric-profile hypothesis, computes the threshold explicitly on geodesic balls of space
forms, and studies its dependence on curvature. The question of Flucher and Rumpf for the
Bernoulli constant $c_K$ itself is not addressed here; for convex $K$, the elliptic solutions of
the interior problem form a family of convex domains \cite{CardaliaguetTahraoui2002,HS00}.

Theorem~3.4 and its corollaries belong to the family of Faber--Krahn type comparison theorems on
Riemannian manifolds. For the first Dirichlet eigenvalue of the Laplacian, the comparison with a
model geodesic ball under a lower Ricci bound is due to Bérard and Meyer \cite{BM82}, and it was
extended to the $p$-Laplacian by Matei \cite{Mat00}; the equality cases and the extension to
non-smooth spaces were studied by Mondino and Semola \cite{MS20}. The threshold considered here
differs from these eigenvalues in one respect: the optimal functions have a plateau $\{u=1\}$ of
positive volume. This affects the rearrangement argument only mildly, but it changes the analysis
of the equality cases, since the rigidity theorem of \cite{MS20} does not apply to functions with
a plateau (Remark~7.7).

The main purpose of this paper is to connect this threshold with the isoperimetric geometry of the ambient space. The first step is a capacity characterization. The second is a P\'olya--Szeg\H{o} comparison based on the isoperimetric profiles. In model spaces, the capacity of concentric geodesic balls is explicitly computable, which leads to a one-dimensional formula for the threshold.

We emphasize that our global comparison theorem is stated under an explicit isoperimetric-profile hypothesis. A lower bound on Ricci curvature is not, by itself, a substitute for this hypothesis on arbitrary noncompact manifolds. We give a counterexample to such an inference and, separately, explain the valid normalized comparison supplied by the L\'evy--Gromov theorem in the compact positive-Ricci setting.

Our main results are the following.
\begin{enumerate}[label=(\roman*)]
\item The threshold is the infimum of a capacity-to-volume quotient (Proposition~2.5), and this
infimum is attained on every bounded Lipschitz domain (Proposition~7.1).
\item Under profile domination, $\Lambda_{p,M}(K)\ge\Lambda_{p,\kappa}(B_R^\kappa)$
(Theorem~3.4); equality forces $K$ to be an isoperimetric region, and a geodesic ball when $M$
is a space form (Theorem~7.5).
\item On geodesic balls of space forms, the threshold is given by a one-dimensional minimization
(Theorem~4.3) and is strictly decreasing in the curvature, at fixed radius and at fixed volume
(Theorems~5.3 and~5.5).
\item On closed manifolds with $\mathrm{Ric}\ge(N-1)\kappa>0$, the normalized Lévy--Gromov
inequality gives a sharp comparison (Corollary~6.1), with equality only for round spheres and
geodesic balls (Theorem~7.6).
\item Under $\mathrm{Ric}\ge0$ and Euclidean volume growth,
$\Lambda_{p,M}(K)\ge\mbox{AVR}^{1/N}\Lambda_{p,\R^N}(K^*)$ (Corollary~6.8), with equality only for
Euclidean space and balls (Theorem~7.8); the cylinder of Example~6.4 shows that the volume growth
assumption cannot be dropped.
\end{enumerate}

%The paper is organized as follows. Section~\ref{sec:variational} gives the variational formulation and proves the capacity characterization. Section~\ref{sec:rearrangement} establishes the rearrangement inequality and the main comparison theorem. Section~\ref{sec:spaceforms} computes the threshold in space forms. Section~\ref{sec:curvature} proves strict curvature monotonicity. Section~\ref{sec:ricci} discusses the positive-Ricci application and the nonnegative-Ricci counterexample. We conclude by stating precisely which perturbative questions are not addressed here.

The paper is organized as follows. Section~2 compares the threshold with the Bernoulli constant
on balls, gives the variational formulation and proves the capacity characterization. Section~3
establishes the Pólya--Szegő inequality under profile domination and the main comparison theorem.
Section~4 computes the threshold on geodesic balls of space forms, and Section~5 proves its strict
monotonicity in the curvature. Section~6 contains the applications to Ricci curvature: the
normalized Lévy--Gromov comparison, the cylinder counterexample, and the comparison under
nonnegative Ricci curvature and Euclidean volume growth. Section~7 proves that the threshold is
attained and characterizes the equality cases. Section~8 collects examples and states precisely
which perturbative questions are not addressed here.

%============================================
\section{Variational formulation and capacity characterization}
\label{sec:variational}
%============================
\subsection{Variational threshold and Bernoulli constant on balls}\label{subsec:ball}

(The formulas used below are derived in Section~4; see \eqref{eq:4-12}--\eqref{eq:4-13}.)

For $0<r<R$, let $u_r$ be the $p$-capacitary potential of $\overline{B_r}$ in $B_R\subset\R^N$
($\Delta_pu_r=0$ in $B_R\setminus\overline{B_r}$, $u_r=1$ on $\partial B_r$, $u_r=0$ on
$\partial B_R$), and let $c(r)=|\nabla u_r|$ on $\partial B_r$. The ball $B_r$ solves the interior
problem with parameter $\lambda$ if and only if $c(r)=\lambda$, so that the Bernoulli constant of
$B_R$ is $c_{B_R}=\min_{0<r<R}c(r)$.

\begin{proposition}\label{prop:ball}
Let $1<p<\infty$ and $R>0$, and let $\rho_*=x_{N,p}R$ be the optimal radius of \eqref{eq:4-13}. Then:
\begin{enumerate}[label=(\roman*),leftmargin=*]
\item $\Lambda_{p,0}(B_R)=c(\rho_*)$; in other words, the critical plateau $B_{\rho_*}$ is a radial
solution of the interior problem with parameter $\Lambda_{p,0}(B_R)$;
\item $c$ attains its minimum at a unique radius $r_*$, with
$r_*=R\big(\tfrac{N-1}{p-1}\big)^{\frac{p-1}{p-N}}$ if $p\neq N$ and $r_*=R/e$ if $p=N$, and
$\rho_*>r_*$; thus $B_{\rho_*}$ lies on the elliptic branch ($c'(\rho_*)>0$);
\item consequently $\Lambda_{p,0}(B_R)>c_{B_R}$, and
\[
\frac{\Lambda_{p,0}(B_R)}{c_{B_R}}=
\begin{cases}
\dfrac{N}{p-1}\Big(\dfrac pN\Big)^{\frac{p-1}{p-N}}\Big(\dfrac{N-1}{p-1}\Big)^{\frac{N-1}{p-N}}, & p\neq N,\\[3mm]
\dfrac{N\,e^{-1/N}}{N-1}, & p=N.
\end{cases}
\]
\end{enumerate}
\end{proposition}

\begin{proof}
Integrating the radial equation gives, for $p\neq N$ and $m=\frac{p-N}{p-1}$,
\[
u_r(\rho)=\frac{R^m-\rho^m}{R^m-r^m},\qquad c(r)=\frac{m\,r^{m-1}}{R^m-r^m},
\]
and $c(r)=1/(r\ln(R/r))$ for $p=N$. In both cases
$\frac{\mathrm d}{\mathrm dr}\ln c(r)=c(r)-\frac{N-1}{(p-1)r}$, so the critical points of $c$ are the
solutions of $c(r)=\frac{N-1}{(p-1)r}$. Since $c(r)\to+\infty$ at both ends of $(0,R)$, a direct
computation shows that this equation has the single solution $r_*$ given in (ii), which is the
minimum of $c$.

(i) For $p\neq N$, \eqref{eq:4-12} gives $x^m=N/p$ for $x=x_{N,p}$, hence
\[
c(xR)=\frac{m\,x^{m-1}}{R(1-x^m)}=\frac{m\,(N/p)\,x^{-1}}{R\,(p-N)/p}=\frac{N}{(p-1)\,x\,R}=\Lambda_{p,0}(B_R).
\]
For $p=N$, $x=e^{-(N-1)/N}$ and $c(xR)=\frac{N}{(N-1)xR}=\Lambda_{N,0}(B_R)$.

(ii) If $p>N$, then $m>0$ and $\rho_*>r_*$ is equivalent to $\frac Np>\frac{N-1}{p-1}$, that is
$N(p-1)>p(N-1)$, that is $p>N$. If $p<N$, then $m<0$ and $\rho_*>r_*$ is equivalent to
$\frac Np<\frac{N-1}{p-1}$, that is $p<N$. If $p=N$, $e^{-(N-1)/N}>e^{-1}$. Since $r_*$ is the only
critical point of $c$ and $\rho_*>r_*$, we have $c'(\rho_*)>0$.

(iii) Since $\rho_*\neq r_*$ and the minimum of $c$ is attained only at $r_*$,
$\Lambda_{p,0}(B_R)=c(\rho_*)>c(r_*)=c_{B_R}$. The formula follows from (i), \eqref{eq:4-13} and
$c_{B_R}=c(r_*)=\frac{N-1}{(p-1)r_*}$.
\end{proof}

\begin{remark}
Proposition~\ref{prop:ball} agrees with Proposition~\ref{prop:critical-attainment}: at $\lambda=\Lambda_{p,0}(B_R)$ the infimum
in \eqref{eq:4-7} is attained, and the corresponding plateau $B_{\rho_*}$ is a nonzero global minimizer, hence
a solution of the free boundary problem. For $c_{B_R}\le\lambda<\Lambda_{p,0}(B_R)$, the interior
problem has radial solutions, but none of them is a global minimizer of \eqref{eq:1-1}. The ratio is
$2e^{-1/2}\approx1.213$ for $N=p=2$ and $1.125$ for $(N,p)=(3,2)$ and $(2,3)$; it decreases with $p$
in Table~\ref{tab:ratio}, and tends to $1$ as $p\to\infty$ since $\rho_*\sim NR/p$ and
$r_*\sim(N-1)R/p$.
\end{remark}
\newpage
\begin{table}[ht]\centering\small
\begin{tabular}{lcccccccc}%\toprule
$(N,p)$ & $(2,1.5)$ & $(2,2)$ & $(2,3)$ & $(2,4)$ & $(3,1.5)$ & $(3,2)$ & $(3,3)$ & $(3,4)$\\%\midrule
$\rho_*/R$ & 0.7500 & 0.6065 & 0.4444 & 0.3536 & 0.7937 & 0.6667 & 0.5134 & 0.4219\\
$r_*/R$ & 0.5000 & 0.3679 & 0.2500 & 0.1925 & 0.6300 & 0.5000 & 0.3679 & 0.2963\\
$\Lambda_{p,0}(B_R)/c_{B_R}$ & 1.3333 & 1.2131 & 1.1250 & 1.0887 & 1.1906 & 1.1250 & 1.0748 & 1.0535\\%\bottomrule
\end{tabular}
\caption{Optimal plateau radius $\rho_*$, turning radius $r_*$ and ratio between the variational
threshold and the Bernoulli constant of $B_R$ (independent numerical check of
Proposition~\ref{prop:ball}).}\label{tab:ratio}
\end{table}

%===========================================================
%\subsection{Variational formulation and capacity characterization in the Riemannian case}
%\label{subsec:variational}
%===========================================================

\subsection{The variational problem in the Riemannian case}
Throughout this section, \(K\subset M\) is a bounded Lipschitz domain in a smooth Riemannian manifold \((M,g)\) of dimension \(N\geq2\), and \(1<p<\infty\). Here $K\Subset M$ means that the closure of $K$ is compact.
%The compact containment is understood to mean that \(\overline K\) is compact.

For \(\lambda>0\), define
\begin{equation}
\mathcal J_{\lambda,p}^{K}(u)
=
\int_K|\nabla_g u|_g^p\,\dd V_g
+(p-1)\lambda^p\Vol_g(\{u<1\}),
\qquad u\in W_0^{1,p}(K).
\label{eq:2-1}\end{equation}
The associated threshold is
\begin{equation}
\Lambda_{p,M}(K)
=
\inf\left\{\lambda>0:
\mathcal J_{\lambda,p}^{K}
\text{ has a nonzero global minimizer in }W_0^{1,p}(K)\right\}.
\label{eq:2-2}\end{equation}
We use the convention that the infimum of the empty set is \(+\infty\). Since
\begin{equation}
\mathcal J_{\lambda,p}^{K}(0)
=(p-1)\lambda^p\Vol_g(K),
\label{eq:2-3}\end{equation}
the threshold measures when a nontrivial competitor can become energetically preferable to zero.

\begin{lemma}[Truncation]
\label{lem:truncation}
For every \(u\in W_0^{1,p}(K)\), the truncation
\[
T(u)=\min\{\max\{u,0\},1\}
\]
belongs to \(W_0^{1,p}(K)\) and satisfies
\[
\mathcal J_{\lambda,p}^{K}(T(u))
\leq \mathcal J_{\lambda,p}^{K}(u).
\]
\end{lemma}

\begin{proof}
The Sobolev chain rule for Lipschitz truncations gives \(T(u)\in W_0^{1,p}(K)\) and
\[
|\nabla_gT(u)|_g\leq|\nabla_gu|_g
\quad\text{a.e. in }K.
\]
Moreover, pointwise,
\[
T(u)<1\quad\Longleftrightarrow\quad u<1.
\]
Thus the volume terms are equal and the Dirichlet term does not increase. This proves the assertion.
\end{proof}

It follows that, when looking for a minimizer with energy no greater than that of zero, we may restrict to
\begin{equation}
0\leq u\leq1\quad\text{a.e. in }K.
\label{eq:2-4}\end{equation}

\subsection{Existence above the threshold}

\begin{lemma}[Existence of a global minimizer]
\label{lem:direct-method}
For every \(\lambda>0\), the functional \(\mathcal J_{\lambda,p}^{K}\) attains its infimum on \(W_0^{1,p}(K)\).
\end{lemma}

\begin{proof}
The functional is bounded below by zero. Let \((u_j)\) be a minimizing sequence. Since
\[
\int_K|\nabla_g u_j|_g^p\,\dd V_g
\leq \mathcal J_{\lambda,p}^{K}(u_j),
\]
the gradients are bounded in \(L^p(K)\). The Poincar\'e inequality on the bounded Lipschitz domain \(K\) then shows that \((u_j)\) is bounded in \(W_0^{1,p}(K)\).

Since \(1<p<\infty\), this space is reflexive. After passing to a subsequence,
\[
u_j\rightharpoonup u\quad\text{weakly in }W_0^{1,p}(K).
\]
The Rellich compactness theorem gives strong convergence in \(L^p(K)\), and, after a further subsequence, \(u_j\to u\) almost everywhere in \(K\).

The Dirichlet integral is weakly lower semicontinuous:
\[
\int_K|\nabla_gu|_g^p\,\dd V_g
\leq
\liminf_{j\to\infty}
\int_K|\nabla_gu_j|_g^p\,\dd V_g.
\]
For the volume term, if \(u(x)<1\), then \(u_j(x)<1\) for all sufficiently large \(j\). Hence
\[
\1_{\{u<1\}}(x)
\leq
\liminf_{j\to\infty}\1_{\{u_j<1\}}(x)
\quad\text{for a.e. }x\in K.
\]
Fatou's lemma yields
\[
\Vol_g(\{u<1\})
\leq
\liminf_{j\to\infty}\Vol_g(\{u_j<1\}).
\]
Combining the two lower-semicontinuity inequalities proves that \(u\) is a global minimizer.
\end{proof}

\subsection{Capacity and the threshold}

For a measurable set \(E\subset K\), define its relative \(p\)-capacity by
\begin{equation}
\RelCap_{p,g}(E,K)
=
\inf\left\{
\int_K|\nabla_gv|_g^p\,\dd V_g:
v\in W_0^{1,p}(K),\
\widetilde v\geq1\quad p\text{-quasi-everywhere on }E
\right\},
\label{eq:2-5}\end{equation}
where \(\widetilde v\) denotes the \(p\)-quasi-continuous representative. We use the standard properties of Sobolev \(p\)-capacity and quasi-continuous representatives; see Maz\'ya~\cite{Mazya}. All quasi-continuity and quasi-everywhere statements are understood relative to \(K\) and to the relative \(p\)-capacity. This definition is meaningful for arbitrary measurable \(E\), since quasi-continuous representatives of a Sobolev function are unique up to a set of relative \(p\)-capacity zero.

Set
\begin{equation}
\mathcal A_p(K)
=
\left\{
u\in W_0^{1,p}(K):
0\leq u\leq1,\ 
\Vol_g(\{u=1\})>0
\right\}.
\label{eq:2-6}\end{equation}

\begin{proposition}[Threshold characterization]
\label{prop:threshold}
The threshold defined in \(\eqref{eq:2-2}\) satisfies
\begin{equation}
\Lambda_{p,M}(K)^p
=
\frac1{p-1}
\inf_{u\in\mathcal A_p(K)}
\frac{\displaystyle\int_K|\nabla_gu|_g^p\,\dd V_g}
{\displaystyle\Vol_g(\{u=1\})}.
\label{eq:2-7}\end{equation}
Moreover,
\begin{equation}
\Lambda_{p,M}(K)^p
=
\frac1{p-1}
\inf_{\substack{E\subset K\ p\text{-quasi-closed}\\ \Vol_g(E)>0}}
\frac{\RelCap_{p,g}(E,K)}{\Vol_g(E)}.
\label{eq:2-8}\end{equation}
\end{proposition}

\begin{proof}
Let
\[
c_K=
\inf_{u\in\mathcal A_p(K)}
\frac{\int_K|\nabla_gu|_g^p\,\dd V_g}
{(p-1)\Vol_g(\{u=1\})}.
\]
For any \(u\in W_0^{1,p}(K)\), Lemma~\ref{lem:truncation} allows us to replace \(u\) by \(T(u)\). For \(0\leq u\leq1\), the partition
\[
K=\{u<1\}\,\dot\cup\,\{u=1\}
\]
holds up to null sets, and therefore
\begin{equation}
\mathcal J_{\lambda,p}^{K}(u)
-\mathcal J_{\lambda,p}^{K}(0)
=
\int_K|\nabla_gu|_g^p\,\dd V_g
-(p-1)\lambda^p\Vol_g(\{u=1\}).
\label{eq:2-9}\end{equation}

Suppose \(0<\lambda^p<c_K\). If \(\Vol_g(\{u=1\})>0\), the definition of \(c_K\) gives
\[
\mathcal J_{\lambda,p}^{K}(u)-\mathcal J_{\lambda,p}^{K}(0)
\geq
(p-1)(c_K-\lambda^p)\Vol_g(\{u=1\})>0.
\]
If \(\Vol_g(\{u=1\})=0\), the difference in \eqref{eq:2-9} is the Dirichlet energy. It vanishes only when \(u\) is constant, and the zero trace then forces \(u=0\). Thus zero is the unique global minimizer below \(c_K^{1/p}\).

Now suppose \(\lambda^p>c_K\). By the definition of the infimum, there exists \(u\in\mathcal A_p(K)\) such that
\[
\int_K|\nabla_gu|_g^p\,\dd V_g
<
(p-1)\lambda^p\Vol_g(\{u=1\}).
\]
Equation \eqref{eq:2-9} gives
\[
\mathcal J_{\lambda,p}^{K}(u)
<
\mathcal J_{\lambda,p}^{K}(0).
\]
By Lemma~\ref{lem:direct-method}, a global minimizer exists. It cannot be zero, because its energy is strictly below the energy of zero. Consequently, every \(\lambda>c_K^{1/p}\) admits a nonzero global minimizer, while no \(\lambda<c_K^{1/p}\) does. This proves \eqref{eq:2-7}. No assertion about attainment at the endpoint \(\lambda=c_K^{1/p}\) is needed.

We prove \eqref{eq:2-8}. Given \(u\in\mathcal A_p(K)\), let
\[
E_u=\{x\in K:\widetilde u(x)=1\},
\]
where \(\widetilde u\) is the \(p\)-quasi-continuous representative. The set \(E_u\) is \(p\)-quasi-closed, since it is the inverse image of the closed set \(\{1\}\) under a \(p\)-quasi-continuous function. Moreover,
$\tilde u=u$ almost everywhere.
%Moreover, \(\widetilde u=u\) outside a set of relative \(p\)-capacity zero, and such a set has Riemannian volume zero.
 Consequently,
\[
\Vol_g(E_u)=\Vol_g(\{u=1\}).
\]
Since \(\widetilde u=1\) on \(E_u\), the function \(u\) is admissible in the definition of \(\RelCap_{p,g}(E_u,K)\). Hence
\[
\RelCap_{p,g}(E_u,K)
\leq\int_K|\nabla_gu|_g^p\,\dd V_g.
\]
Taking infima gives
\begin{equation}
\inf_E\frac{\RelCap_{p,g}(E,K)}{\Vol_g(E)}
\leq
\inf_{u\in\mathcal A_p(K)}
\frac{\int_K|\nabla_gu|_g^p\,\dd V_g}
{\Vol_g(\{u=1\})}.
\label{eq:2-10}\end{equation}

Conversely, let \(E\subset K\) be \(p\)-quasi-closed with positive measure. If
\(\RelCap_{p,g}(E,K)=+\infty\), the desired inequality for this set is
automatic. We may therefore assume that its relative capacity is finite.
For every \(\delta>0\), choose \(v_\delta\in W_0^{1,p}(K)\) admissible
for \(E\) such that
\[
\int_K|\nabla_gv_\delta|_g^p\,\dd V_g
\leq\RelCap_{p,g}(E,K)+\delta.
\]
Truncate \(v_\delta\) to \([0,1]\). The truncated function \(u_\delta\) still equals \(1\) quasi-everywhere on \(E\), so
\[
\Vol_g(\{u_\delta=1\})\geq\Vol_g(E),
\qquad
\int_K|\nabla_gu_\delta|_g^p\,\dd V_g
\leq\RelCap_{p,g}(E,K)+\delta.
\]
Thus
\[
\inf_{u\in\mathcal A_p(K)}
\frac{\int_K|\nabla_gu|_g^p\,\dd V_g}
{\Vol_g(\{u=1\})}
\leq
\frac{\RelCap_{p,g}(E,K)+\delta}{\Vol_g(E)}.
\]
Letting \(\delta\downarrow0\) and taking the infimum over \(E\) gives the reverse inequality to \eqref{eq:2-10}. This proves \eqref{eq:2-8}.
\end{proof}

\begin{corollary}[Nondegeneracy of the threshold]
\label{cor:threshold-positive-finite}
Under the hypotheses of this section,
\[
0<\Lambda_{p,M}(K)<+\infty.
\]
\end{corollary}

\begin{proof}
By the Poincar\'e inequality there is a constant $C_P>0$ such that
\[
\int_K|\nabla_g u|_g^p\,\dd V_g
\ge C_P\int_K|u|^p\,\dd V_g
\ge C_P\Vol_g(\{u=1\})
\]
for every $u\in\mathcal A_p(K)$. Formula \eqref{eq:2-7} yields
$\Lambda_{p,M}(K)^p\ge C_P/(p-1)>0$.
Since $K$ is a nonempty open domain, it contains a relatively compact
geodesic ball. Choose $u\in C_c^\infty(K)$ with $0\le u\le1$ and
$u=1$ on a smaller ball. Then $u\in\mathcal A_p(K)$ and its
energy-to-plateau-volume quotient is finite. Formula \eqref{eq:2-7} gives
$\Lambda_{p,M}(K)<+\infty$.
\end{proof}

\begin{proposition}[Attainment at the critical threshold]
\label{prop:critical-attainment}
Let
\[
c_K=(p-1)\Lambda_{p,M}(K)^p
=
\inf_{u\in\mathcal A_p(K)}
\frac{\displaystyle\int_K|\nabla_g u|_g^p\,\dd V_g}
{\Vol_g(\{u=1\})}.
\]
Then \(u=0\) is a global minimizer of
\(\mathcal J_{\Lambda_{p,M}(K),p}^{K}\). Moreover, the following
statements are equivalent:
\begin{enumerate}
\item \(\mathcal J_{\Lambda_{p,M}(K),p}^{K}\) has a nonzero global
minimizer;
\item the infimum defining \(c_K\) is attained by some
\(u_*\in\mathcal A_p(K)\).
\end{enumerate}
In that case, every such \(u_*\), after truncation to \([0,1]\), is a
global minimizer at the critical parameter and has the same energy as
the zero function.

\end{proposition}

\begin{proof}
By the truncation lemma, it is enough to consider competitors with
\(0\leq u\leq1\). For every such \(u\), identity \eqref{eq:2-9} and the
definition of \(c_K\) give
\[
\mathcal J_{\Lambda_{p,M}(K),p}^{K}(u)
-\mathcal J_{\Lambda_{p,M}(K),p}^{K}(0)
=
\int_K|\nabla_g u|_g^p\,\dd V_g
-c_K\Vol_g(\{u=1\})\geq0.
\]
Thus zero is a global minimizer at the critical parameter.

Suppose first that the infimum defining \(c_K\) is attained by
\(u_*\in\mathcal A_p(K)\). The same identity gives
\[
\mathcal J_{\Lambda_{p,M}(K),p}^{K}(u_*)
=
\mathcal J_{\Lambda_{p,M}(K),p}^{K}(0).
\]
Since zero is already a global minimizer, \(u_*\) is also a global
minimizer, and it is nonzero because
\(\Vol_g(\{u_*=1\})>0\).

Conversely, suppose that \(w\ne0\) is a global minimizer at the
critical parameter. Replace it by \(u=T(w)\). Then \(u\) is also a
global minimizer and \(0\leq u\leq1\). We first note that \(u\ne0\):
otherwise \(w\leq0\) almost everywhere, so that
\(\Vol_g(\{w<1\})=\Vol_g(K)\) and
\[
\mathcal J_{\Lambda_{p,M}(K),p}^{K}(w)
=\int_K|\nabla_gw|_g^p\,\dd V_g
 +(p-1)\Lambda_{p,M}(K)^p\Vol_g(K)
\geq \mathcal J_{\Lambda_{p,M}(K),p}^{K}(0).
\]
Since $w$ is a global minimizer and zero is an admissible competitor,
equality must hold. Equality is possible only if \(\nabla_gw=0\) almost everywhere, which by
Poincar\'e's inequality implies \(w=0\), contrary to the assumption.
Thus \(u\ne0\). Its energy equals that of zero. If
\(\Vol_g(\{u=1\})=0\), identity \eqref{eq:2-9} implies
\(\int_K|\nabla_g u|_g^p\,\dd V_g=0\), and Poincar\'e's inequality
forces \(u=0\), again a contradiction. Hence \(u\in\mathcal A_p(K)\).
Identity \eqref{eq:2-9} now yields
\[
\frac{\displaystyle\int_K|\nabla_g u|_g^p\,\dd V_g}
{\Vol_g(\{u=1\})}
=c_K.
\]
The infimum is therefore attained.
\end{proof}

%\begin{remark}
%This proposition gives an exact criterion for endpoint attainment;
%it does not assert that the capacity-to-volume infimum is attained on
%every bounded Lipschitz domain.  We invite the reader to see  Proposition \ref{prop:attain} for a positive response.

%\end{remark}
\begin{remark}%{Remark 2.8}
Proposition~2.7 gives an exact criterion for endpoint attainment. It does not by itself show
that the capacity-to-volume infimum is attained; this is proved for every bounded Lipschitz
domain in Proposition~7.1.
\end{remark}

%===========================================================
\section{Isoperimetric rearrangement and the main comparison theorem}
\label{sec:rearrangement}
%===========================================================

\subsection{Isoperimetric profiles}

For a Riemannian manifold \((M,g)\), define its isoperimetric profile by
\begin{equation}
I_M(v)
=
\inf\left\{
\Per_g(E):
E\subset M,\ \Vol_g(E)=v
\right\},
\label{eq:3-1}\end{equation}
where \(\Per_g(E)\) is the perimeter in the sense of finite-perimeter sets.

Let \((M_\kappa^N,g_\kappa)\) be the simply connected space form of dimension \(N\) and constant sectional curvature \(\kappa\). Let \(K\Subset M\) be a bounded Lipschitz domain and let \(B_R^\kappa\subset M_\kappa^N\) satisfy
\begin{equation}
\Vol_g(K)=\Vol_{g_\kappa}(B_R^\kappa)=V.
\label{eq:3-2}\end{equation}
We assume the profile domination
\begin{equation}
I_M(v)\geq I_\kappa(v),
\qquad 0<v\leq V,
\label{eq:3-3}\end{equation}
%and that, for every \(v\in(0,V]\), a geodesic ball in the model space realizes \(I_\kappa(v)\). We also assume that these model balls can be chosen as concentric balls inside \(B_R^\kappa\).
%For the simply connected space forms considered in Section~\ref{sec:spaceforms},
%this realization property is the classical isoperimetric theorem in
%constant-curvature spaces; see Chavel~\cite[Chapter~VI]{Chavel}.

By the isoperimetric theorem in simply connected space forms \cite{Chavel,BZ88b}, for every
$v\in(0,V]$ the value $I_\kappa(v)$ is the area of the boundary of a geodesic ball of volume $v$,
and these balls can be chosen concentric inside $B_R^\kappa$. We refer to this as the
\emph{model-ball realization property}; in the present setting it is a theorem, not an
assumption.

\subsection{P\'olya--Szeg\H{o} inequality}

\begin{lemma}[Extension across the pole]
\label{lem:radial-extension}
Let $N\geq2$, let $B_R^\kappa$ be a proper geodesic ball in an
$N$-dimensional space form (in particular, $R<\pi/\sqrt{\kappa}$
when $\kappa>0$), and let
$P_\kappa(r)=\Per_{g_\kappa}(\partial B_r^\kappa)$.
Suppose that $f:(0,R)\to\mathbb R$ is locally absolutely continuous,
bounded, and satisfies
\[
\int_0^R\bigl(|f(r)|^p+|f'(r)|^p\bigr)P_\kappa(r)\,\dd r<\infty.
\]
Then the radial function $F(x)=f(d(o,x))$, initially defined away from
the pole $o$, belongs to $W^{1,p}(B_R^\kappa)$ and
\[
\int_{B_R^\kappa}|\nabla F|^p\,\dd V_{g_\kappa}
=\int_0^R|f'(r)|^pP_\kappa(r)\,\dd r.
\]
If, in addition, $f$ vanishes on a neighborhood of $R$, then
$F\in W_0^{1,p}(B_R^\kappa)$.
\end{lemma}

\begin{proof}
On every annulus $B_R^\kappa\setminus\overline{B_\varepsilon^\kappa}$,
polar coordinates and local absolute continuity give the weak gradient
$\nabla F=f'(r)\partial_r$. The assumed weighted integrability shows
that $F$ and this candidate gradient belong to $L^p$ on the punctured
ball, with the stated norm identity. It remains to check that no
additional distributional derivative is concentrated at the pole.
For a smooth compactly supported test vector field $X$, integration by
parts on the punctured ball produces an inner-boundary term bounded by
\[
\|f\|_{L^\infty(0,R)}\,\|X\|_{L^\infty}\,
P_\kappa(\varepsilon).
\]
Since $P_\kappa(\varepsilon)=O(\varepsilon^{N-1})$ and $N\geq2$,
this term tends to zero as $\varepsilon\downarrow0$. Passing to the
limit proves that the candidate gradient is the distributional
Riemannian gradient on the whole ball. Thus $F\in W^{1,p}(B_R^\kappa)$
and the energy identity follows by polar integration. If $f$ vanishes
near $R$, then $F$ vanishes near the boundary of the ball; the zero
trace characterization of $W_0^{1,p}$ on a Lipschitz ball gives the
last assertion.
\end{proof}

For \(u\geq0\), define its distribution function by
\begin{equation}
\mu_u(t)=\Vol_g(\{u>t\}),\qquad t\geq0.
\label{eq:3-4}\end{equation}
Its spherical decreasing rearrangement \(u^*\) on \(B_R^\kappa\) is the radial nonincreasing function characterized by
\begin{equation}
\Vol_{g_\kappa}(\{u^*>t\})=\mu_u(t)
\quad\text{for every }t\geq0.
\label{eq:3-5}\end{equation}

\begin{lemma}[P\'olya--Szeg\H{o} under profile domination]
\label{lem:PS}
%Assume (3.2)--(3.3) and the model-ball realization hypothesis.
Assume \eqref{eq:3-2} and \eqref{eq:3-3}. For every nonnegative $u\in W^{1,p}_0(K)$, its
spherical decreasing rearrangement $u^*$ belongs to $W^{1,p}_0(B_R^\kappa)$ and satisfies
\begin{equation}
\int_{B_R^\kappa}|\nabla_{g_\kappa}u^*|^p\,\dd V_{g_\kappa}
\leq
\int_K|\nabla_gu|_g^p\,\dd V_g.
\label{eq:3-6}\end{equation}
\end{lemma}

\begin{proof}
We use  \cite{Ta76} and the general rearrangement theorem of Nobili and Violo
\cite[Theorem~1.1]{NobiliViolo}. We first identify its one-dimensional
model space in the present setting.

Let
\[
P_\kappa(r)=\Per_{g_\kappa}(\partial B_r^\kappa),
\qquad 0<r<R,
\]
and equip the interval $(0,R)$ with the measure
\[
\dd\omega_\kappa(r)=P_\kappa(r)\,\dd r.
\]
The polar-coordinate formula gives
\[
\omega_\kappa((0,r))=\int_0^rP_\kappa(s)\,\dd s
=\Vol_{g_\kappa}(B_r^\kappa),
\qquad \omega_\kappa((0,R))=V.
\]
In particular, the cumulative distribution function of $\omega_\kappa$
is strictly increasing, and its inverse sends $v\in(0,V)$ to the unique
radius $r_v\in(0,R)$ such that
$\Vol_{g_\kappa}(B_{r_v}^\kappa)=v$. The profile associated with this
weighted interval, in the notation of \cite[Section~3.2]{NobiliViolo},
is therefore
\[
I_{\omega_\kappa}^{\flat}(v)=P_\kappa(r_v)=I_\kappa(v),
\qquad 0<v<V,
\]
where the last equality is exactly the model-ball realization
property.
%model-ball realization hypothesis.

We first take $u\in C_c^\infty(K)$, $u\geq0$. Choose an open
neighborhood $U$ of $\operatorname{supp}u$ with
$\overline U\subset K$ and $\Vol_g(U)<V$. Such a choice is possible
because $\operatorname{supp}u$ is compactly contained in $K$. For every
finite-perimeter set $E\subset U$,
its perimeter relative to $K$ agrees with its ambient perimeter in $M$.
Consequently, \eqref{eq:3-3} and the preceding identification imply
\[
\Per_g(E;K)\geq I_M(\Vol_g(E))
\geq I_{\omega_\kappa}^{\flat}(\Vol_g(E)).
\]
Since $\operatorname{supp}u$ is compactly contained in $U$, we may
choose $\varepsilon>0$ so small that
$(\operatorname{supp}u)^\varepsilon\cap K\subset U$.
Thus the perimeter hypothesis in \cite[Theorem~1.1]{NobiliViolo}
holds on the required $\varepsilon$-neighborhood of the support,
with constant equal to one. Applying that theorem, with the target
weighted interval $((0,R),\omega_\kappa)$, yields a decreasing
rearrangement $u^\sharp$ satisfying
\[
\int_0^R |(u^\sharp)'(r)|^pP_\kappa(r)\,\dd r
\leq \int_K|\nabla_gu|_g^p\,\dd V_g.
\]
By construction, $u^\sharp(r)$ is the radial profile of the spherical
decreasing rearrangement $u^*$; polar coordinates identify the left-hand
side with
$\int_{B_R^\kappa}|\nabla_{g_\kappa}u^*|^p\,\dd V_{g_\kappa}$.
Moreover, since $u$ is compactly supported in $K$, its positive support
has volume strictly smaller than $V$. Hence $u^*$ vanishes on an outer
annulus of $B_R^\kappa$. The rearrangement is bounded by
$\|u\|_{L^\infty(K)}$ and has finite weighted radial energy. Lemma
\ref{lem:radial-extension}, applied to its radial profile, therefore
gives $u^*\in W^{1,p}(B_R^\kappa)$ with the asserted energy identity.
Its vanishing near the outer boundary then gives
$u^*\in W_0^{1,p}(B_R^\kappa)$. This proves \eqref{eq:3-6} for nonnegative
smooth compactly supported functions.

Now let $u\in W_0^{1,p}(K)$ be nonnegative. Choose nonnegative
$u_j\in C_c^\infty(K)$ such that
$u_j\to u$ strongly in $W_0^{1,p}(K)$. The inequality just proved
shows that $(u_j^*)$ is bounded in $W_0^{1,p}(B_R^\kappa)$. The
contraction property of decreasing rearrangement on spaces of equal
finite measure gives
\[
\|u_j^*-u^*\|_{L^p(B_R^\kappa)}
\leq \|u_j-u\|_{L^p(K)}\longrightarrow0.
\]
After extraction of a subsequence, $u_j^*$ converges weakly in
$W_0^{1,p}(B_R^\kappa)$; the strong $L^p$ convergence identifies its
weak limit with $u^*$. By weak lower semicontinuity,
\[
\begin{aligned}
\int_{B_R^\kappa}|\nabla_{g_\kappa}u^*|^p\,\dd V_{g_\kappa}
&\leq \liminf_{j\to\infty}
\int_{B_R^\kappa}|\nabla_{g_\kappa}u_j^*|^p\,\dd V_{g_\kappa}\\
&\leq \lim_{j\to\infty}
\int_K|\nabla_gu_j|_g^p\,\dd V_g
=\int_K|\nabla_gu|_g^p\,\dd V_g.
\end{aligned}
\]
This proves the lemma.
\end{proof}

\begin{remark}
%The model-ball realization assumption
The model-ball realization property  is essential for the stated radial
rearrangement. A profile inequality without a model geometry in which
the profile is realized by nested balls does not, by itself, produce the
particular rearrangement used here. The proof above uses the general
Sobolev rearrangement theorem of Nobili--Violo rather than a generic
Morse perturbation argument; no assumption on nondegeneracy or distinct
critical values is needed.
\end{remark}

\subsection{Main comparison theorem}

\begin{theorem}[Isoperimetric comparison of variational thresholds]
%\begin{remark}
%Before proving this theorem, let us note that, in a simply
%connected space form, isoperimetric regions are geodesic balls \cite{Chavel}, and the balls of
%volume $v\le V$ can be chosen concentric inside $B^\kappa_R.$\\
%In the case of equality, the standard tool of rigidity can be found in  \cite{BZ88}, see also\cite{Chavel}.
%\end{remark}
\label{thm:main}
%Let $(M,g)$ be a complete $N$-dimensional Riemannian manifold, $K\Subset M$ a bounded Lipschitz
%domain, $1<p<\infty$, $\kappa\in\mathbb R$, and let $B_R^\kappa\subset M_\kappa^N$ satisfy (3.2),
%with $R<\pi/\sqrt\kappa$ if $\kappa>0$. If the profile domination (3.3) holds, then

Let $(M,g)$ be a complete $N$-dimensional Riemannian manifold, $K\Subset M$ a bounded Lipschitz
domain, $1<p<\infty$ and $\kappa\in\R$, and let $B_R^\kappa\subset M_\kappa^N$ satisfy
\eqref{eq:3-2}, with $R<\pi/\sqrt\kappa$ if $\kappa>0$. If the profile domination \eqref{eq:3-3}
holds, then
\begin{equation}
\Lambda_{p,M}(K)\ge\Lambda_{p,\kappa}(B_R^\kappa). 
\label{eq:3-10}\end{equation}
\end{theorem}

\begin{proof}
By Proposition~\ref{prop:threshold},
\begin{equation}
\Lambda_{p,M}(K)^p
=
\frac1{p-1}
\inf_{u\in\mathcal A_p(K)}
\frac{\int_K|\nabla_gu|_g^p\,\dd V_g}
{\Vol_g(\{u=1\})}.
\label{eq:3-11}\end{equation}
Fix \(u\in\mathcal A_p(K)\). By Lemma~\ref{lem:PS}, its spherical decreasing rearrangement \(u^*\) belongs to \(W_0^{1,p}(B_R^\kappa)\) and satisfies
\begin{equation}
\int_{B_R^\kappa}|\nabla_{g_\kappa}u^*|^p\,\dd V_{g_\kappa}
\leq
\int_K|\nabla_gu|_g^p\,\dd V_g.
\label{eq:3-12}\end{equation}
Equimeasurability gives
\[
\Vol_{g_\kappa}(\{u^*>t\})
=
\Vol_g(\{u>t\})
\quad (t\geq0).
\]
Since \(0\leq u,u^*\leq1\), continuity of measure from above, applied as \(t\uparrow1\), yields
\begin{equation}
\Vol_{g_\kappa}(\{u^*=1\})
=
\Vol_g(\{u=1\})>0.
\label{eq:3-13}\end{equation}
Thus \(u^*\in\mathcal A_p(B_R^\kappa)\). By the threshold characterization on the model ball,
\[
\begin{aligned}
\Lambda_{p,\kappa}(B_R^\kappa)^p
&\leq
\frac1{p-1}
\frac{\displaystyle\int_{B_R^\kappa}
|\nabla_{g_\kappa}u^*|^p\,\dd V_{g_\kappa}}
{\Vol_{g_\kappa}(\{u^*=1\})}\\
&\leq
\frac1{p-1}
\frac{\displaystyle\int_K|\nabla_gu|_g^p\,\dd V_g}
{\Vol_g(\{u=1\})}.
\end{aligned}
\]
Taking the infimum over \(u\in\mathcal A_p(K)\) and using \eqref{eq:3-11} proves \eqref{eq:3-10}.
\end{proof}
\begin{remark}
Equality holds in \eqref{eq:3-10} when $M=M_\kappa^N$ and $K=B_R^\kappa$ (Section~8.1).
Conversely, Theorem~7.5 shows that equality forces $K$ to be an isoperimetric region of $M$ with
$\Per_g(K)=I_\kappa(\Vol_g(K))$, and a geodesic ball when $M=M_\kappa^N$. The proof relies on the
analysis of the equality case in the Pólya--Szegő inequality, in the spirit of Brothers and
Ziemer \cite{BZ88}.
\end{remark}

%===========================================================
\section{Exact computation on geodesic balls in space forms}
\label{sec:spaceforms}
%===========================================================

\subsection{Polar coordinates and volume}

Let \(M_\kappa^N\) be the simply connected space form of dimension \(N\) and constant sectional curvature \(\kappa\). Define
\begin{equation}
S_\kappa(r)=
\begin{cases}
\dfrac{\sin(\sqrt\kappa\,r)}{\sqrt\kappa},&\kappa>0,\\[2mm]
r,&\kappa=0,\\[1mm]
\dfrac{\sinh(\sqrt{-\kappa}\,r)}{\sqrt{-\kappa}},&\kappa<0.
\end{cases}
\label{eq:4-1}\end{equation}
In geodesic polar coordinates,
\begin{equation}
g_\kappa=\dd r^2+S_\kappa(r)^2g_{\mathbb S^{N-1}},
\qquad
\dd V_{g_\kappa}
=S_\kappa(r)^{N-1}\dd r\,\dd\sigma.
\label{eq:4-2}\end{equation}
We write
\begin{equation}
V_\kappa(R)
=
\omega_{N-1}\int_0^R S_\kappa(t)^{N-1}\,\dd t.
\label{eq:4-3}\end{equation}
For \(\kappa>0\), the radius is restricted to \(0<R<\pi/\sqrt\kappa\).

\subsection{Relative capacity of concentric balls}

\begin{proposition}
\label{prop:capacity-model}
For \(0<\rho<R\),
\begin{equation}
\RelCap_{p,\kappa}(B_\rho^\kappa,B_R^\kappa)
=
\omega_{N-1}
\left(
\int_\rho^R
S_\kappa(t)^{-\frac{N-1}{p-1}}\,\dd t
\right)^{1-p}.
\label{eq:4-4}\end{equation}
\end{proposition}

\begin{proof}
We first justify the reduction to radial functions. Let \(v\) be
admissible for the capacity and replace it by its truncation to
\([0,1]\). Then \(v=1\) almost everywhere in \(B_\rho^\kappa\).
For almost every \(r\in(\rho,R)\), define the spherical average
\[
\overline v(r)
=
\frac{1}{\omega_{N-1}}
\int_{\mathbb S^{N-1}}v(r,\theta)\,\dd\theta.
\]
The Sobolev slicing and trace theorems imply that
\(\overline v\in W^{1,p}((\rho,R))\), that its endpoint traces satisfy
\[
\overline v(\rho)=1,\qquad \overline v(R)=0,
\]
and that
\[
\overline v'(r)
=
\frac{1}{\omega_{N-1}}
\int_{\mathbb S^{N-1}}\partial_rv(r,\theta)\,\dd\theta
\quad\text{for a.e. }r.
\]
Jensen's inequality therefore gives
\[
|\overline v'(r)|^p
\leq
\frac{1}{\omega_{N-1}}
\int_{\mathbb S^{N-1}}|\partial_rv(r,\theta)|^p\,\dd\theta
\leq
\frac{1}{\omega_{N-1}}
\int_{\mathbb S^{N-1}}|\nabla v(r,\theta)|^p\,\dd\theta.
\]
The radial function equal to \(1\) in \(B_\rho^\kappa\) and to
\(\overline v(r)\) on the annulus is consequently admissible and has
no greater energy. Thus it suffices to minimize over radial functions.

For any absolutely continuous radial competitor \(w\) with
\(w(\rho)=1\) and \(w(R)=0\), Hölder's inequality yields
\[
\begin{aligned}
1
&=\left|\int_\rho^R w'(r)\,\dd r\right|\\
&\leq
\left(\int_\rho^R |w'(r)|^pS_\kappa(r)^{N-1}\,\dd r\right)^{1/p}
\left(\int_\rho^R
S_\kappa(r)^{-\frac{N-1}{p-1}}\,\dd r
\right)^{(p-1)/p}.
\end{aligned}
\]
It follows that
\[
\int_\rho^R |w'(r)|^pS_\kappa(r)^{N-1}\,\dd r
\geq
\left(\int_\rho^R
S_\kappa(r)^{-\frac{N-1}{p-1}}\,\dd r
\right)^{1-p}.
\]
Equality is attained by the absolutely continuous function
\[
w_*(r)=
\frac{\displaystyle\int_r^R
S_\kappa(t)^{-\frac{N-1}{p-1}}\,\dd t}
{\displaystyle\int_\rho^R
S_\kappa(t)^{-\frac{N-1}{p-1}}\,\dd t},
\qquad \rho\leq r\leq R,
\]
extended by \(1\) on \(B_\rho^\kappa\). Multiplication by
\(\omega_{N-1}\) gives \eqref{eq:4-4}.
\end{proof}

\subsection{Threshold formula}

\begin{lemma}[Balls minimize relative capacity at fixed volume in a model ball]
\label{lem:capacity-rearrangement}
Let \(E\subset B_R^\kappa\) be \(p\)-quasi-closed with
\[
0<\Vol_{g_\kappa}(E)=\Vol_{g_\kappa}(B_\rho^\kappa)
\]
for some \(0<\rho<R\). Then
\[
\RelCap_{p,\kappa}(B_\rho^\kappa,B_R^\kappa)
\leq
\RelCap_{p,\kappa}(E,B_R^\kappa).
\]
\end{lemma}

\begin{proof}
Let \(v\in W_0^{1,p}(B_R^\kappa)\) be admissible for
\(\RelCap_{p,\kappa}(E,B_R^\kappa)\). Replacing \(v\) by its truncation to
\([0,1]\) does not increase its energy, so we may assume \(0\leq v\leq1\)
and \(v=1\) quasi-everywhere on \(E\). A set of relative \(p\)-capacity
zero has Riemannian volume zero. Thus, for every \(t<1\),
\[
\Vol_{g_\kappa}(\{v>t\})\geq\Vol_{g_\kappa}(E)
=\Vol_{g_\kappa}(B_\rho^\kappa).
\]
Let \(v^*\) be the spherical decreasing rearrangement of \(v\) in
\(B_R^\kappa\). By equimeasurability and the fact that \(0\leq v^*\leq1\), the
preceding inequality for every \(t<1\) implies
\(v^*=1\) almost everywhere on \(B_\rho^\kappa\). A Sobolev function
which equals \(1\) almost everywhere on an open set has
\(p\)-quasi-continuous representative equal to \(1\)
\(p\)-quasi-everywhere there (apply the corresponding locality
property to \(v^*-1\)). Hence \(v^*\) is admissible for the
capacity of \(B_\rho^\kappa\). %The
%P\'olya--Szeg\H{o} inequality (see \cite{Ta76}) gives
The Pólya--Szegő inequality in $M_\kappa^N$ (Lemma~3.2 with $M=M_\kappa^N$ and $K=B_R^\kappa$,
for which \eqref{eq:3-3} holds with equality) gives
\[
\RelCap_{p,\kappa}(B_\rho^\kappa,B_R^\kappa)
\leq\int_{B_R^\kappa}|\nabla v^*|^p\,\dd V_{g_\kappa}
\leq\int_{B_R^\kappa}|\nabla v|^p\,\dd V_{g_\kappa}.
\]
Taking the infimum over all admissible \(v\) proves the lemma.
\end{proof}

\begin{theorem}[Exact threshold on a model ball]
\label{thm:threshold-ball}
For \(1<p<\infty\),
\begin{equation}
\Lambda_{p,\kappa}(B_R^\kappa)^p
=
\inf_{0<\rho<R}
\frac{\omega_{N-1}}
{(p-1)V_\kappa(\rho)}
\left(
\int_\rho^R
S_\kappa(t)^{-\frac{N-1}{p-1}}\,\dd t
\right)^{1-p}.
\label{eq:4-7}\end{equation}
The infimum is attained at an interior radius \(\rho_*\in(0,R)\).
\end{theorem}

\begin{proof}
By Proposition~\ref{prop:threshold}, the threshold is the infimum of the relative capacity-to-volume quotient over \(p\)-quasi-closed positivity sets of positive volume. Lemma~\ref{lem:capacity-rearrangement} shows that, among sets of prescribed volume, a concentric geodesic ball has no larger relative \(p\)-capacity. Hence the infimum reduces to \(B_\rho^\kappa\), \(0<\rho<R\). Proposition~\ref{prop:capacity-model} and \eqref{eq:4-3} give \eqref{eq:4-7}.

We prove attainment. The function of \(\rho\) in \eqref{eq:4-7} is continuous on \((0,R)\). As \(\rho\uparrow R\), the integral from \(\rho\) to \(R\) tends to zero; since \(1-p<0\), the quotient tends to \(+\infty\).

As \(\rho\downarrow0\), we use \(S_\kappa(\rho)\sim\rho\). Put \(q=(N-1)/(p-1)\). If \(p<N\), then \(q>1\) and
\[
\int_\rho^R S_\kappa(t)^{-q}\,\dd t
\asymp \rho^{1-q},
\]
so the capacity is of order \(\rho^{N-p}\), whereas \(V_\kappa(\rho)\asymp\rho^N\). Their quotient is of order \(\rho^{-p}\), which tends to \(+\infty\). If \(p=N\), the integral is of order \(|\log\rho|\), and the quotient is of order \(\rho^{-N}|\log\rho|^{1-N}\), which also tends to \(+\infty\). If \(p>N\), the integral has a finite positive limit as \(\rho\downarrow0\), while the volume is of order \(\rho^N\); the quotient again tends to \(+\infty\).

The quotient therefore tends to \(+\infty\) at both endpoints and is continuous in between. It attains its minimum at some \(\rho_*\in(0,R)\).
\end{proof}

\subsection{Dimensionless representation}

Let
\[
\tau=\kappa R^2,\qquad x=\frac{\rho}{R},
\]
and define
\begin{equation}
s_\tau(t)=
\begin{cases}
\dfrac{\sin(\sqrt\tau\,t)}{\sqrt\tau},&\tau>0,\\[2mm]
t,&\tau=0,\\[1mm]
\dfrac{\sinh(\sqrt{-\tau}\,t)}{\sqrt{-\tau}},&\tau<0.
\end{cases}
\label{eq:4-8}\end{equation}
For positive curvature, admissibility requires \(0<\sqrt\tau<\pi\). Set
\begin{equation}
A_\tau(x)=\int_0^x s_\tau(t)^{N-1}\,\dd t,
\qquad
B_\tau(x)=\int_x^1s_\tau(t)^{-q}\,\dd t,
\qquad q=\frac{N-1}{p-1},
\label{eq:4-9}\end{equation}
and
\begin{equation}
F_\tau(x)
=
\frac{1}{(p-1)A_\tau(x)}B_\tau(x)^{1-p},
\qquad
m(\tau)=\min_{0<x<1}F_\tau(x).
\label{eq:4-10}\end{equation}
Changing variables \(t=Rs\) in \eqref{eq:4-7} gives
\begin{equation}
\Lambda_{p,\kappa}(B_R^\kappa)^p
=
R^{-p}m(\kappa R^2).
\label{eq:4-11}\end{equation}

\subsection{The Euclidean model and an explicit example}

When \(\kappa=0\), one has \(s_0(t)=t\). The minimizing radius \(x_{N,p}\) is found by differentiating
\[
F_0(x)
=
\frac{1}{(p-1)\displaystyle\int_0^x t^{N-1}\,\dd t}
\left(\int_x^1t^{-q}\,\dd t\right)^{1-p}.
\]
Since \(\int_0^x t^{N-1}\dd t=x^N/N\), the critical-point equation is
\[
-\frac{N}{x}
+\frac{p-1}{x^q\displaystyle\int_x^1t^{-q}\,\dd t}=0.
\]
Equivalently,
\begin{equation}
\int_x^1t^{-q}\,\dd t
=
\frac{p-1}{N}x^{1-q}.
\label{eq:4-12}\end{equation}
If \(p\ne N\), evaluating the integral gives
\[
x^{1-q}=\frac Np,
\qquad
x_{N,p}
=
\left(\frac Np\right)^{\frac{p-1}{p-N}}.
\]
If \(p=N\), then \(q=1\), and \eqref{eq:4-12} becomes
\[
-\log x=\frac{N-1}{N},
\qquad
x_{N,N}=e^{-\frac{N-1}{N}}.
\]
Substituting \eqref{eq:4-12} into the quotient yields
\begin{equation}
\Lambda_{p,0}(B_R)
=
\frac{N}{(p-1)R\,x_{N,p}},
\qquad
x_{N,p}=
\begin{cases}
\left(\dfrac Np\right)^{\frac{p-1}{p-N}},&p\ne N,\\[2mm]
e^{-\frac{N-1}{N}},&p=N.
\end{cases}
\label{eq:4-13}\end{equation}

\begin{example}[The case \(N=3\), \(p=2\)]
Here
\[
x_{3,2}=\left(\frac32\right)^{-1}=\frac23.
\]
Thus the optimal inner radius is \(\rho_*=2R/3\), and
\begin{equation}
\Lambda_{2,0}(B_R)=\frac{9}{2R}.
\label{eq:4-14}\end{equation}
This provides a concrete benchmark for the curvature-dependent formula \eqref{eq:4-7}.
\end{example}

%===========================================================
\section{Strict monotonicity with respect to curvature}
\label{sec:curvature}
%===========================================================

We prove that the model threshold decreases strictly when the sectional curvature increases. We first keep the geodesic radius fixed and then keep the volume fixed.

\subsection{A differential identity}

For admissible \(\tau\) and \(t>0\), define
\begin{equation}
a_\tau(t)=-\partial_\tau\log s_\tau(t).
\label{eq:5-1}\end{equation}

\begin{lemma}
\label{lem:a}
For every admissible \(\tau\),
\begin{equation}
a_\tau(t)>0,
\qquad
a_\tau'(t)
=
\frac{\displaystyle\int_0^t s_\tau(u)^2\,\dd u}
{s_\tau(t)^2}>0.
\label{eq:5-2}\end{equation}
\end{lemma}

\begin{proof}
The function \(s_\tau\) satisfies
\[
s_\tau''+\tau s_\tau=0,\qquad
s_\tau(0)=0,\qquad s_\tau'(0)=1.
\]
Let \(c=s_\tau'/s_\tau\). Then
\[
c'+c^2+\tau=0.
\]
Writing \(b=\partial_\tau c\) and differentiating gives
\[
b'+2cb=-1.
\]
Since \(a_\tau=-\partial_\tau\log s_\tau\), we have \(a_\tau'=-b\). Multiplying the equation for \(b\) by \(s_\tau^2\) gives
\[
(b s_\tau^2)'=-s_\tau^2.
\]
The Taylor expansion of the solution of
\(s_\tau''+\tau s_\tau=0\), with its initial conditions, is
\[
s_\tau(t)=t-\frac{\tau}{6}t^3+O(t^5),
\qquad
\partial_\tau s_\tau(t)=-\frac16t^3+O(t^5).
\]
Consequently,
\[
c(t)=\frac{s_\tau'(t)}{s_\tau(t)}
=\frac1t-\frac{\tau}{3}t+O(t^3),
\qquad
b(t)=\partial_\tau c(t)=-\frac13t+O(t^3).
\]
In particular,
\[
b(t)s_\tau(t)^2=O(t^3)\longrightarrow0
\quad\text{as }t\downarrow0.
\]
Integration from \(0\) to \(t\) therefore yields
\[
b(t)s_\tau(t)^2=-\int_0^t s_\tau(u)^2\,\dd u.
\]
Thus
\[
a_\tau'(t)=-b(t)
=
\frac{\int_0^t s_\tau(u)^2\,\dd u}{s_\tau(t)^2}>0.
\]
The same expansion gives
\[
a_\tau(t)
=-\frac{\partial_\tau s_\tau(t)}{s_\tau(t)}
=\frac16t^2+O(t^4),
\]
so \(a_\tau(t)\to0\) as \(t\downarrow0\). Since \(a_\tau'(t)>0\),
we conclude that \(a_\tau(t)>0\) for every \(t>0\).
\end{proof}

\subsection{Monotonicity at fixed radius}

For fixed \(x\in(0,1)\), define
\[
\alpha_\tau(x)
=
\frac{\displaystyle\int_0^x
a_\tau(t)s_\tau(t)^{N-1}\,\dd t}
{\displaystyle\int_0^x s_\tau(t)^{N-1}\,\dd t},
\]
and
\[
\beta_\tau(x)
=
\frac{\displaystyle\int_x^1
a_\tau(t)s_\tau(t)^{-q}\,\dd t}
{\displaystyle\int_x^1s_\tau(t)^{-q}\,\dd t}.
\]
Since \(a_\tau\) is strictly increasing,
\begin{equation}
\alpha_\tau(x)<a_\tau(x)<\beta_\tau(x).
\label{eq:5-3}\end{equation}

\begin{lemma}
\label{lem:Fderivative}
For every \(x\in(0,1)\),
\begin{equation}
\partial_\tau\log F_\tau(x)
=
(N-1)\big(\alpha_\tau(x)-\beta_\tau(x)\big)<0.
\label{eq:5-4}\end{equation}
\end{lemma}

\begin{proof}
From \(a_\tau=-\partial_\tau\log s_\tau\), differentiation under the integral sign gives
\[
\partial_\tau A_\tau(x)
=-(N-1)\int_0^x
a_\tau(t)s_\tau(t)^{N-1}\,\dd t.
\]
Hence
\[
\partial_\tau\log A_\tau(x)=-(N-1)\alpha_\tau(x).
\]
Similarly,
\[
\partial_\tau B_\tau(x)
=q\int_x^1a_\tau(t)s_\tau(t)^{-q}\,\dd t,
\]
so
\[
\partial_\tau\log B_\tau(x)=q\beta_\tau(x).
\]
Using
\[
\log F_\tau=-\log(p-1)-\log A_\tau+(1-p)\log B_\tau
\]
and \((p-1)q=N-1\), we obtain
\[
\partial_\tau\log F_\tau
=(N-1)\alpha_\tau-(N-1)\beta_\tau<0.
\]
\end{proof}

\begin{theorem}[Strict curvature monotonicity at fixed radius]
\label{thm:fixed-radius}
Fix \(R>0\). If \(\kappa_1<\kappa_2\) and both geodesic balls of radius \(R\) are admissible in their respective space forms, then
\begin{equation}
\Lambda_{p,\kappa_2}(B_R^{\kappa_2})
<
\Lambda_{p,\kappa_1}(B_R^{\kappa_1}).
\label{eq:5-5}\end{equation}
\end{theorem}

\begin{proof}
Set
\[
\tau_i=\kappa_iR^2,\qquad i=1,2.
\]
Since \(R>0\) and \(\kappa_1<\kappa_2\), we have
\[
\tau_1<\tau_2.
\]
The admissibility assumption means that the model warping functions
\(s_{\tau_i}\) are positive on \((0,1)\). In the positively curved
case this is equivalent to \(\tau_i<\pi^2\). Consequently, every
\(\tau\in[\tau_1,\tau_2]\) is also admissible: if \(\tau_2\leq0\)
there is no positive-curvature restriction, while if \(\tau_2>0\)
then \(\tau\leq\tau_2<\pi^2\). Thus \(F_\tau(x)\) is well-defined
throughout this parameter interval for every \(x\in(0,1)\).

We first record the pointwise comparison. Fix \(x\in(0,1)\).
Lemma~\ref{lem:Fderivative} applies for every
\(\tau\in[\tau_1,\tau_2]\) and gives
\[
\frac{\partial}{\partial\tau}\log F_\tau(x)
=
(N-1)\bigl(\alpha_\tau(x)-\beta_\tau(x)\bigr)<0.
\]
The map \(\tau\mapsto F_\tau(x)\) is continuously differentiable
on this interval. Integrating the preceding strict inequality, we
obtain
\[
\log F_{\tau_2}(x)-\log F_{\tau_1}(x)
=
\int_{\tau_1}^{\tau_2}
(N-1)\bigl(\alpha_\tau(x)-\beta_\tau(x)\bigr)\,\dd\tau
<0.
\]
Since the exponential function is strictly increasing, this proves
\begin{equation}
F_{\tau_2}(x)<F_{\tau_1}(x)
\qquad\text{for every }x\in(0,1).
\label{eq:5-6}\end{equation}

We now pass from the pointwise comparison to the minimum. By
Theorem~\ref{thm:threshold-ball}, applied to the two admissible
model balls, the function \(F_{\tau_1}\) attains its minimum at
some point \(x_1\in(0,1)\). In particular,
\[
m(\tau_1)=F_{\tau_1}(x_1).
\]
No uniqueness of \(x_1\) is needed; we fix any one of its minimizers.
By the definition of \(m(\tau_2)\) as the infimum (in fact, the
minimum) of \(F_{\tau_2}\) over \((0,1)\), and by \eqref{eq:5-6} evaluated
at this fixed point \(x_1\), we have
\[
\begin{aligned}
m(\tau_2)
&=\min_{0<x<1}F_{\tau_2}(x)\\
&\leq F_{\tau_2}(x_1)\\
&<F_{\tau_1}(x_1)\\
&=m(\tau_1).
\end{aligned}
\]
Hence
\begin{equation}
m(\tau_2)<m(\tau_1).
\label{eq:5-7}\end{equation}
This argument does not differentiate \(m\) and does not require
the minimizing point to depend continuously or differentiably on
\(\tau\). This is important because uniqueness or smooth dependence
of the minimizing radius has not been established in the paper.

Finally, the dimensionless identity \eqref{eq:4-11} gives, since \(R\) is
the same for the two model balls,
\[
\Lambda_{p,\kappa_i}(B_R^{\kappa_i})^p
=R^{-p}m(\kappa_iR^2)
=R^{-p}m(\tau_i),\qquad i=1,2.
\]
Combining this identity with \eqref{eq:5-7}, and using \(R^{-p}>0\), yields
\[
\Lambda_{p,\kappa_2}(B_R^{\kappa_2})^p
<
\Lambda_{p,\kappa_1}(B_R^{\kappa_1})^p.
\]
Both thresholds are nonnegative, so taking the positive \(p\)-th
root preserves the strict inequality. This proves \eqref{eq:5-5}.
\end{proof}

\subsection{Monotonicity at fixed volume}

Fix \(V>0\). For each admissible curvature \(\kappa\), let \(R_\kappa\) be the radius of a model ball of volume \(V\). The admissible curvatures form an interval: all \(\kappa\leq0\) are allowed, while for \(\kappa>0\) we require \(V<\Vol(S_\kappa^N)\). For each such \(\kappa\), the radius \(R_\kappa\) is unique because \(R\mapsto V_\kappa(R)\) is strictly increasing. Put
\[
\tau=\kappa R_\kappa^2,\qquad
H(\tau)=\int_0^1s_\tau(t)^{N-1}\,\dd t.
\]
The volume constraint is
\begin{equation}
V=\omega_{N-1}R_\kappa^N H(\tau).
\label{eq:5-8}\end{equation}

\begin{lemma}
\label{lem:fixed-volume-derivatives}
Along a differentiable fixed-volume family,
\begin{equation}
\frac{\dd\tau}{\dd\kappa}
=
\frac{N R_\kappa^2}
{N+2\tau H'(\tau)/H(\tau)}>0,
\qquad
%\[
\frac{R_\kappa'}{R_\kappa}
=
-\frac{R_\kappa^2 H'(\tau)/H(\tau)}
{N+2\tau H'(\tau)/H(\tau)}>0.
\label{eq:5-9}\end{equation}
\end{lemma}

\begin{proof}
The function
\[
\mathcal V(\kappa,R)
=\omega_{N-1}\int_0^R S_\kappa(t)^{N-1}\,\dd t
\]
is continuously differentiable (indeed smooth) in the admissible
parameter range, including at \(\kappa=0\), as follows from the
smooth dependence of \(S_\kappa\) on \(\kappa\). Moreover,
\[
\partial_R\mathcal V(\kappa,R)
=\omega_{N-1}S_\kappa(R)^{N-1}>0.
\]
The implicit function theorem therefore gives a locally smooth
function \(R_\kappa\) satisfying
\(\mathcal V(\kappa,R_\kappa)=V\). Uniqueness of the radius makes
these local functions agree on overlaps, so \(R_\kappa\) is
differentiable throughout the common admissible curvature interval.

Lemma~\ref{lem:a} gives
\begin{equation}
H'(\tau)
=-(N-1)\int_0^1a_\tau(t)s_\tau(t)^{N-1}\,\dd t<0.
\label{eq:5-10}\end{equation}
Write \(y=R_\kappa'/R_\kappa\). Differentiating \(\tau=\kappa R_\kappa^2\) gives
\begin{equation}
\tau'=R_\kappa^2+2\tau y.
\label{eq:5-11}\end{equation}
Differentiating the fixed-volume identity \eqref{eq:5-8} and dividing by \(V\) gives
\begin{equation}
0=Ny+\frac{H'(\tau)}{H(\tau)}\tau'.
\label{eq:5-12}\end{equation}
Solving \eqref{eq:5-11}--\eqref{eq:5-12}, we find
\[
\tau'
=
\frac{N R_\kappa^2}{N+2\tau H'(\tau)/H(\tau)}
\]
and
\[
y
=
-\frac{R_\kappa^2H'(\tau)/H(\tau)}
{N+2\tau H'(\tau)/H(\tau)}.
\]

It remains to verify that the common denominator is positive. At fixed \(\kappa\), the volume of a ball is strictly increasing with its radius:
\[
\frac{\partial}{\partial R}V_\kappa(R)
=
\omega_{N-1}S_\kappa(R)^{N-1}>0.
\]
On the other hand, differentiating the expression
\[
V_\kappa(R)=\omega_{N-1}R^NH(\kappa R^2)
\]
with respect to \(R\) gives
\[
\frac{\partial}{\partial R}V_\kappa(R)
=
\omega_{N-1}R^{N-1}
\big(NH(\tau)+2\tau H'(\tau)\big).
\]
Therefore
\[
N+2\tau H'(\tau)/H(\tau)>0.
\]
Since \(H'(\tau)<0\), both expressions in \eqref{eq:5-9} are strictly positive.
\end{proof}

\begin{theorem}[Strict curvature monotonicity at fixed volume]
\label{thm:fixed-volume}
Fix \(V>0\) and consider a common admissible interval of curvatures for which the geodesic balls of volume \(V\) exist. If \(\kappa_1<\kappa_2\), then
\begin{equation}
\Lambda_{p,\kappa_2}(B_{R_{\kappa_2}}^{\kappa_2})
<
\Lambda_{p,\kappa_1}(B_{R_{\kappa_1}}^{\kappa_1}).
\label{eq:5-13}\end{equation}
\end{theorem}

\begin{proof}
Let \(\kappa_1<\kappa_2\) belong to the common admissible
curvature interval, and write
\[
R_i=R_{\kappa_i},\qquad
\tau_i=\kappa_iR_i^2,\qquad i=1,2.
\]
We shall compare separately the geometric scale \(R_\kappa^{-1}\)
and the dimensionless factor \(m(\tau(\kappa))^{1/p}\).

\medskip
\noindent
\textit{Step 1: strict increase of the radius.}
By Lemma~\ref{lem:fixed-volume-derivatives}, the fixed-volume
radius is differentiable and satisfies
\[
\frac{R_\kappa'}{R_\kappa}
=
-\frac{R_\kappa^2H'(\tau)/H(\tau)}
{N+2\tau H'(\tau)/H(\tau)}.
\]
The denominator is positive by the last part of that lemma.
Moreover, \(H(\tau)>0\) and \(H'(\tau)<0\), the latter following
from \eqref{eq:5-10} and the strict positivity of \(a_\tau\) on \((0,1)\).
Since \(R_\kappa>0\), it follows that
\[
R_\kappa'>0
\]
throughout the admissible interval. Integrating this strict
derivative inequality on \([\kappa_1,\kappa_2]\), we obtain
\begin{equation}
0<R_1<R_2,
\qquad\text{and hence}\qquad
0<\frac{R_1}{R_2}<1.
\label{eq:5-14}\end{equation}

\medskip
\noindent
\textit{Step 2: strict increase of the dimensionless curvature
parameter.}
The same lemma gives
\[
\tau'(\kappa)
=
\frac{N R_\kappa^2}
{N+2\tau H'(\tau)/H(\tau)}>0.
\]
Consequently,
\begin{equation}
\tau_1<\tau_2.
\label{eq:5-15}\end{equation}
This conclusion is not obtained merely from the increase of
\(\kappa\), since the radius also depends on \(\kappa\); it is the
fixed-volume differential identity \eqref{eq:5-9} that establishes it.

\medskip
\noindent
\textit{Step 3: strict decrease of the dimensionless minimum.}
Theorem~\ref{thm:fixed-radius} proves that \(m\) is strictly
decreasing on the admissible dimensionless curvature interval.
%Applying that result to (5.14), we get
Applying that result to \eqref{eq:5-15}, we get
%$0<m(\tau_2)<m(\tau_1)$.
\[
0<m(\tau_2)<m(\tau_1).
\]
The positivity follows from the positivity of \(F_\tau(x)\) and
the existence of an interior minimizer established in
Section~4. Therefore
\begin{equation}
0<
\left(\frac{m(\tau_2)}{m(\tau_1)}\right)^{1/p}
<1.
\label{eq:5-16}\end{equation}

\medskip
\noindent
\textit{Step 4: comparison of the thresholds.}
The dimensionless identity \eqref{eq:4-11}, applied to the two
fixed-volume balls, gives
\[
\Lambda_{p,\kappa_i}(B_{R_i}^{\kappa_i})
=
R_i^{-1}m(\tau_i)^{1/p},
\qquad i=1,2.
\]
%Taking the ratio and using (5.12) and (5.14), we find
Taking the ratio and using \eqref{eq:5-14} and \eqref{eq:5-16}, we find
\[
\begin{aligned}
\frac{\Lambda_{p,\kappa_2}(B_{R_2}^{\kappa_2})}
{\Lambda_{p,\kappa_1}(B_{R_1}^{\kappa_1})}
&=
\frac{R_1}{R_2}
\left(\frac{m(\tau_2)}{m(\tau_1)}\right)^{1/p}\\
&<1.
\end{aligned}
\]
This proves
\[
\Lambda_{p,\kappa_2}(B_{R_{\kappa_2}}^{\kappa_2})
<
\Lambda_{p,\kappa_1}(B_{R_{\kappa_1}}^{\kappa_1}),
\]
as claimed.

Notice that no derivative of \(m\) has been used. In particular,
the proof does not require uniqueness of the minimizing
dimensionless radius, nor differentiability of that minimizer
with respect to \(\tau\). The strict inequality follows from the
two monotonicity statements already established: the radius
increases, while the dimensionless minimum decreases.
\end{proof}

\begin{remark}
For positive curvature and fixed volume, the admissible interval is restricted by the requirement that the prescribed volume be smaller than the total volume of the sphere of curvature \(\kappa\). The fixed-volume theorem is asserted only on a common interval where the model balls exist.
\end{remark}

%===========================================================
\section{Ricci curvature: a valid positive result and a counterexample}
\label{sec:ricci}
%===========================================================

\subsection{The normalized L\'evy--Gromov comparison}

The profile hypothesis in Theorem~\ref{thm:main} can be verified in important geometric situations. One such situation is the compact positive-Ricci setting. The correct formulation uses normalized volume measures.

Let \((M^N,g)\) be closed and connected, and suppose
\begin{equation}
\Ric_g\geq(N-1)\kappa g,\qquad \kappa>0.
\label{eq:6-1}\end{equation}
Let \((S_\kappa^N,g_\kappa)\) be the round sphere of sectional curvature \(\kappa\). Write
\[
V_M=\Vol_g(M),\qquad V_\kappa^{\mathrm{tot}}=\Vol_{g_\kappa}(S_\kappa^N).
\]
Since \(S_\kappa^N\) is the round sphere of sectional curvature
\(\kappa>0\), its polar-coordinate expression gives
\[
V_\kappa^{\mathrm{tot}}
=
\omega_{N-1}\int_0^{\pi/\sqrt{\kappa}}
S_\kappa(t)^{N-1}\,\dd t,
\qquad
S_\kappa(t)=\frac{\sin(\sqrt{\kappa}\,t)}{\sqrt{\kappa}},
\]
where \(\omega_{N-1}=\Vol(S^{N-1}(1))\) is the area of the
unit \((N-1)\)-dimensional sphere. Equivalently,
\[
V_\kappa^{\mathrm{tot}}
=\kappa^{-N/2}\Vol(S^N(1)).
\]
We then define the probability measures
\[
\mu_M=\frac{\dd V_g}{V_M},
\qquad
\mu_\kappa=\frac{\dd V_{g_\kappa}}{V_\kappa^{\mathrm{tot}}}.
\]
For \(0<\theta<1\), let
\[
\mathcal I_M(\theta)
=
\inf\left\{
\frac{\Per_g(E)}{V_M}:
\mu_M(E)=\theta
\right\},
\]
and define \(\mathcal I_\kappa\) analogously on the sphere.

The L\'evy--Gromov isoperimetric inequality states that
\begin{equation}
\mathcal I_M(\theta)\geq\mathcal I_\kappa(\theta),
\qquad 0<\theta<1.
\label{eq:6-2}\end{equation}
This is a comparison of normalized profiles, not an assertion that the unnormalized profile of \(M\) dominates the profile of a volume-matched round sphere. See \cite{CavallettiMondino,Bayle, Gro99}.

\begin{corollary}[Threshold comparison under positive Ricci curvature]
\label{cor:levy-gromov}
Assume \eqref{eq:6-1}, and let \(K\Subset M\) be a bounded Lipschitz domain with
\[
0<\theta:=\frac{\Vol_g(K)}{\Vol_g(M)}<1.
\]
Let \(B_\theta^\kappa\subset S_\kappa^N\) be a geodesic ball satisfying
\[
\frac{\Vol_{g_\kappa}(B_\theta^\kappa)}
{\Vol_{g_\kappa}(S_\kappa^N)}
=\theta.
\]
If \(r_\theta\in(0,\pi/\sqrt{\kappa})\) is chosen by
\[
\frac{\displaystyle\int_0^{r_\theta}S_\kappa(t)^{N-1}\,\dd t}
{\displaystyle\int_0^{\pi/\sqrt{\kappa}}S_\kappa(t)^{N-1}\,\dd t}
=\theta,
\]
then \(B_\theta^\kappa=B_{r_\theta}^\kappa\). The comparison reads
\begin{equation}
\Lambda_{p,M}(K)\geq
\Lambda_{p,\kappa}(B_{r_\theta}^\kappa).
\label{eq:6-3}\end{equation}
\end{corollary}

\begin{proof}
We apply the rearrangement theorem of Nobili--Violo
\cite[Theorem~1.1]{NobiliViolo} with probability measures. Put
$L=\pi/\sqrt{\kappa}$ and define on $(0,L)$ the weighted measure
\[
\dd\omega_\kappa(r)
=\frac{P_\kappa(r)}{V_\kappa^{\mathrm{tot}}}\,\dd r,
\qquad
P_\kappa(r)=\omega_{N-1}S_\kappa(r)^{N-1}.
\]
Its total mass is one, and its cumulative distribution is
\[
\omega_\kappa((0,r))
=\frac{\Vol_{g_\kappa}(B_r^\kappa)}
{V_\kappa^{\mathrm{tot}}}.
\]
Consequently, the inverse cumulative distribution sends
$\theta\in(0,1)$ to $r_\theta$, and the associated weighted-interval
profile is
\[
I_{\omega_\kappa}^{\flat}(\theta)
=\frac{P_\kappa(r_\theta)}{V_\kappa^{\mathrm{tot}}}
=\mathcal I_\kappa(\theta),
\]
where the last equality is the isoperimetric profile of the round
sphere, realized by geodesic balls (or their complements).

First take $u\in C_c^\infty(K)$, $u\geq0$. Choose an open set
$U\Subset K$ containing $\operatorname{supp}u$. For every
finite-perimeter set $E\subset U$, its perimeter relative to $K$
coincides with its ambient perimeter in $M$. The normalized
L\'evy--Gromov inequality therefore gives
\[
\Per_g(E;K)/V_M
\geq \mathcal I_M(\mu_M(E))
\geq \mathcal I_\kappa(\mu_M(E))
=I_{\omega_\kappa}^{\flat}(\mu_M(E)).
\]
Since $\operatorname{supp}u$ is compactly contained in $U$, choose
$\varepsilon>0$ such that
$(\operatorname{supp}u)^\varepsilon\cap K\subset U$.
Thus the perimeter hypothesis of
\cite[Theorem~1.1]{NobiliViolo} holds on the required
$\varepsilon$-neighborhood of $\operatorname{supp}u$, with constant
one, for the metric measure space $(M,d_g,\mu_M)$ and target weighted interval
$((0,L),\omega_\kappa)$. Since $\mu_M(K)=\theta$, the canonical
rearrangement interval has measure $\theta$ and is $(0,r_\theta)$.
The theorem yields a decreasing rearrangement $u^\star$ such that
\[
\int_0^{r_\theta}|(u^\star)'(r)|^p\,
\frac{P_\kappa(r)}{V_\kappa^{\mathrm{tot}}}\,\dd r
\leq
\int_K|\nabla_gu|_g^p\,\dd\mu_M.
\]
Under polar coordinates, $u^\star$ is the radial rearrangement on
$B_{r_\theta}^\kappa$, and the left-hand side is its normalized
Dirichlet energy. Because $u$ is compactly supported in $K$, its
rearrangement vanishes near the outer boundary of this ball. 
%Its
%boundedness and finite weighted radial energy give the Sobolev
%extension across the center allow us to apply
Since it is
bounded and has finite weighted radial energy, Lemma~\ref{lem:radial-extension} applies and gives the Sobolev extension
across the centre; hence
%Lemma~\ref{lem:radial-extension}; hence
$u^\star\in W_0^{1,p}(B_{r_\theta}^\kappa)$.

For general nonnegative $u\in W_0^{1,p}(K)$, choose nonnegative
$u_j\in C_c^\infty(K)$ converging strongly to $u$ in $W_0^{1,p}(K)$.
The preceding inequality bounds the rearrangements in
$W_0^{1,p}(B_{r_\theta}^\kappa)$. The $L^p$ contraction of decreasing
rearrangement for probability spaces gives strong $L^p$ convergence
to the rearrangement $u^\star$ of $u$. Weak compactness and weak lower
semicontinuity therefore yield
\[
\int_{B_{r_\theta}^\kappa}
|\nabla_{g_\kappa}u^\star|^p\,\dd\mu_\kappa
\leq
\int_K|\nabla_gu|^p\,\dd\mu_M.
\]
Equimeasurability gives
$\mu_\kappa(\{u^\star=1\})=\mu_M(\{u=1\})$ for
$0\leq u\leq1$. Dividing the energy and the volume by their
respective total volumes does not change their quotient:
\[
\frac{\int_K|\nabla_gu|^p\,\dd\mu_M}
{\mu_M(\{u=1\})}
=
\frac{\int_K|\nabla_gu|^p\,\dd V_g}
{\Vol_g(\{u=1\})},
\]
and similarly on the model ball. Thus the normalized quotient cannot
increase under rearrangement. Taking the infimum in the threshold
characterization proves \eqref{eq:6-3}.
\end{proof}

\begin{remark}
The model ball in Corollary~\ref{cor:levy-gromov} has the same \emph{volume fraction} as \(K\), not necessarily the same unnormalized volume. This distinction is necessary because the L\'evy--Gromov theorem compares normalized isoperimetric profiles.
\end{remark}

\begin{example}[The round sphere]
If \(M=S_\kappa^N\) and \(K=B_R^\kappa\) is a geodesic ball, the normalized L\'evy--Gromov comparison is an equality. Consequently, the threshold comparison in Corollary~\ref{cor:levy-gromov} is an equality for the corresponding volume fraction. Its value is given by the explicit formula in Theorem~\ref{thm:threshold-ball}.
\end{example}

\subsection{Why \texorpdfstring{\(\Ric\geq0\)}{Ricci nonnegative} alone does not imply a global Euclidean profile comparison}

The following example explains why the nonnegative-Ricci assumption cannot simply replace the profile hypothesis \eqref{eq:3-3} with the Euclidean profile on arbitrary noncompact manifolds.

\begin{example}[A product-cylinder counterexample]
\label{ex:cylinder}
Let \(N\geq2\) and consider
\[
M=S^{N-1}(a)\times\R
\]
with the product metric, where \(S^{N-1}(a)\) is the round sphere of radius \(a>0\). The Ricci tensor of the product is nonnegative: in the sphere directions it equals \((N-2)a^{-2}g\), and in the \(\R\)-direction it vanishes. Thus
\[
\Ric_M\geq0.
\]

Let
\[
K_L=S^{N-1}(a)\times(-L,L).
\]
Writing \(A=\Vol(S^{N-1}(a))\), we have
\begin{equation}
\Vol_M(K_L)=2LA,
\qquad
\Per_M(K_L)=2A.
\label{eq:6-6}\end{equation}
Let \(B_L\subset\R^N\) be the Euclidean ball with the same volume. Its radius is
\[
r_L=\left(\frac{2LA}{\omega_N}\right)^{1/N},
\]
and its perimeter is
\begin{equation}
\Per_{\R^N}(B_L)
=
N\omega_N r_L^{N-1}
=
N\omega_N^{1/N}(2LA)^{(N-1)/N}.
\label{eq:6-7}\end{equation}
Consequently,
\begin{equation}
\frac{\Per_M(K_L)}{\Per_{\R^N}(B_L)}
=
\frac{2A}
{N\omega_N^{1/N}(2LA)^{(N-1)/N}}
\longrightarrow0
\quad\text{as }L\to+\infty.
\label{eq:6-8}\end{equation}
For all sufficiently large \(L\),
\[
\Per_M(K_L)<\Per_{\R^N}(B_L).
\]
Therefore the Euclidean isoperimetric lower bound fails on this manifold, despite \(\Ric_M\geq0\).
\end{example}

\begin{remark}
Example~\ref{ex:cylinder} disproves the implication
\[
\Ric_g\geq0
\quad\Longrightarrow\quad
I_M(v)\geq I_{\R^N}(v)\ \text{for all }v>0
\]
in the noncompact setting. It does not, by itself, disprove a possible variational-threshold inequality under some other hypotheses. It shows that the proof of such an inequality cannot be based solely on the asserted global Euclidean profile domination.
\end{remark}

\subsection{Nonnegative Ricci curvature and Euclidean volume growth}\label{sec:avr}

%This section can be inserted as a new Section~6.3, after Example~6.4.
 Let $(M,g)$ be a complete,
noncompact, $N$-dimensional Riemannian manifold with $\mathrm{Ric}_g\ge0$. By the Bishop--Gromov
theorem, the ratio $\Vol_g(B_r(o))/(\beta_Nr^N)$ is nonincreasing in $r$, and its limit
\[
%\AVR_g 
\mbox{AVR}_g = \lim_{r\to\infty}\frac{\Vol_g(B_r(o))}{\beta_N\,r^N}\in[0,1]
\]
does not depend on $o$. One says that $M$ has \emph{Euclidean volume growth} if $\mbox{AVR}_g>0.$%$\AVR_g>0$. 
The
cylinder of Example~6.4 has linear volume growth, hence $\mbox{AVR}_g=0$.

\begin{lemma}[Isoperimetric input]\label{lem:iso}
Assume $\mathrm{Ric}_g\ge0$ and $\theta:=\mbox{AVR}_g>0$. For every set $E\subset M$ of finite perimeter
with compact closure,
\[
\Per_g(E)\ \ge\ \theta^{1/N}\,I_0\big(\Vol_g(E)\big),\qquad I_0(v)=N\beta_N^{1/N}v^{\frac{N-1}N}.
\]
\end{lemma}

\begin{proof}
This is the sharp isoperimetric inequality of Balogh and Kristály \cite{BK23}, stated there for
bounded Borel sets with the outer Minkowski content in place of the perimeter (the Riemannian case
was also proved by Brendle \cite{Bre21}). For a set of finite perimeter with compact closure, it follows by
approximating $E$ by smooth sets $E_k$ with $\Vol_g(E_k\triangle E)\to0$ and
$\Per_g(E_k)\to\Per_g(E)$, for which the Minkowski content and the perimeter coincide.
\end{proof}

\begin{lemma}[Pólya--Szegő inequality with a constant]\label{lem:ps}
Under the assumptions of Lemma~\ref{lem:iso}, let $K\Subset M$ be a bounded Lipschitz domain and let
$B\subset\R^N$ be the Euclidean ball centred at the origin with $|B|=\Vol_g(K)$. For every nonnegative
$u\in W^{1,p}_0(K)$, its Schwarz rearrangement $u^\sharp$ on $B$ belongs to $W^{1,p}_0(B)$ and
%\begin{equation}\label{eq:ps}
\begin{equation}\theta^{p/N}\int_B|\nabla u^\sharp|^p\dd x\ \le\ \int_K|\nabla_gu|_g^p\dd V_g .\label{eq:6-9}\end{equation}
%\end{equation}
\end{lemma}

\begin{proof}
The proof is that of Lemma~3.2, with two changes. First, the target weighted interval is
$\big((0,R_B),N\beta_Nr^{N-1}\dd r\big)$, whose isoperimetric profile is $I_0$. Second, by
Lemma~6.6, the perimeter hypothesis of \cite[Theorem~1.1]{NobiliViolo} holds on every open set
$U\Subset K$ with the multiplicative constant $C=\theta^{1/N}$. Theorem~1.1 of \cite{NobiliViolo}
allows such a constant, and its conclusion then carries the factor $C^p=\theta^{p/N}$. This
gives \eqref{eq:6-9} for $u\in C^\infty_c(K)$, $u\ge0$. The extension to $W^{1,p}_0(K)$ is the
approximation argument at the end of the proof of Lemma~3.2.
\end{proof}

\begin{corollary}[Threshold comparison under $\mathrm{Ric}\ge0$]\label{cor:avr}
Let $(M,g)$ be a complete, noncompact, $N$-dimensional Riemannian manifold with
$\mathrm{Ric}_g\ge0$ and Euclidean volume growth, let $K\Subset M$ be a bounded Lipschitz domain, and
let $B\subset\R^N$ be a Euclidean ball with $|B|=\Vol_g(K)$, of radius
$R_K^{\rm e}=(\Vol_g(K)/\beta_N)^{1/N}$. Then
\[
\Lambda_{p,M}(K)\ \ge\ \mbox{AVR}_g^{1/N}\,\Lambda_{p,0}(B)
=\mbox{AVR}_g^{1/N}\,\frac{N}{(p-1)\,x_{N,p}\,R_K^{\rm e}} .
\]
\end{corollary}

\begin{proof}
Let $\theta=\mbox{AVR}_g$ and $u\in A_p(K)$. By Lemma~\ref{lem:ps}, $u^\sharp\in W^{1,p}_0(B)$, and, as
in \eqref{eq:3-13}, $|\{u^\sharp=1\}|=\Vol_g(\{u=1\})>0$, so that $u^\sharp\in A_p(B)$. By
Proposition \ref{prop:threshold} on $B$ and \eqref{eq:6-9},
\[
(p-1)\,\Lambda_{p,0}(B)^p\ \le\ \frac{\int_B|\nabla u^\sharp|^p\dd x}{|\{u^\sharp=1\}|}
\ \le\ \theta^{-p/N}\,\frac{\int_K|\nabla_gu|_g^p\dd V_g}{\Vol_g(\{u=1\})} .
\]
Taking the infimum over $u\in A_p(K)$ and using \eqref{eq:2-7} gives
$\Lambda_{p,0}(B)^p\le\theta^{-p/N}\Lambda_{p,M}(K)^p$, which is the claim.
\end{proof}

\begin{remark}[Comparison with Example~6.4 and with Theorem~3.4]
Corollary~\ref{cor:avr} is the positive counterpart of Example~6.4: the cylinder has $\mbox{AVR}=0$, and
the corollary gives no information there. When $M=\R^N$, $\mbox{AVR}=1$ and the corollary reduces to the
Euclidean case of Theorem~3.4. For $\mbox{AVR}_g<1$, the hypothesis \eqref{eq:3-3} of Theorem~3.4 with $\kappa=0$
fails for large volumes, which is why a constant factor is unavoidable.
\end{remark}

\begin{remark}[On the optimality of the factor $\mbox{AVR}^{1/N}$]
%\begin{paperrem}[On the optimality of the factor $\AVR^{1/N}$]{Remark 6.10}
Let $\Sigma=S^{N-1}(a)$ with $0<a<1$, and let $C(\Sigma)$ be the metric cone over $\Sigma$, a
$\mathrm{CD}(0,N)$ space with $\mbox{AVR}=a^{N-1}=:\theta$, singular at its vertex $o$. Radial
competitors on $B_r(o)$ have exactly $\theta$ times the energy and $\theta$ times the plateau
volume of their Euclidean counterparts on $B_r$, so the computation of Section~4 gives
$\Lambda_p(B_r(o))\le\Lambda_{p,0}(B_r)=N/((p-1)x_{N,p}r)$. On the other hand,
$\Vol(B_r(o))=\theta\beta_Nr^N$, so the right-hand side of Corollary~6.8 for $K=B_r(o)$ is
$\theta^{1/N}\cdot N/((p-1)x_{N,p}\theta^{1/N}r)=N/((p-1)x_{N,p}r)$. Hence, if Corollary~6.8
extends to such cones (its two inputs, \cite{BK23} and \cite{NobiliViolo}, are available on metric
measure spaces), equality holds for balls centred at the vertex, and the factor $\mbox{AVR}^{1/N}$
cannot be improved in general. On smooth manifolds, equality occurs only for Euclidean space and
balls (Theorem~7.8).
%\end{paperrem}
\end{remark}

%===============================
\section{Attainment of the threshold and equality cases}\label{sec:rig}
%==============================================

%This section can be inserted as a new Section~7, before the examples.
 We first show that the
infimum in \eqref{eq:2-7} is attained on every bounded Lipschitz domain, as announced in Remark~2.8. We then characterize the equality cases in Theorem~3.4 when $M$ is a space
form, in Corollary~6.1, and in Corollary~6.8. The comparison of Corollary~6.1 is the analogue, for the variational threshold,
of the Bérard--Meyer comparison for the first Dirichlet eigenvalue \cite{BM82}, whose equality
case (for the $p$-Laplacian, and on non-smooth spaces) is treated by Mondino and Semola
\cite[Section~5]{MS20}. For $u\in A_p(K)$ we write
\[
\mathcal R_K(u)=\frac{\int_K|\nabla_gu|_g^p\dd V_g}{\Vol_g(\{u=1\})},
\qquad\text{so that}\qquad (p-1)\,\Lambda_{p,M}(K)^p=\inf_{u\in A_p(K)}\mathcal R_K(u)
\]
by Proposition \ref{prop:threshold}; the same notation is used on model balls. In the setting of Corollary~6.1,
volumes are normalized by $\Vol_g(M)$ and $\Vol_{g_\kappa}(S^N_\kappa)$; since $\mathcal R$ is a
quotient, this changes none of the statements below.

\subsection{Attainment}

\begin{proposition}[The threshold is attained]\label{prop:attain}
Let $(M,g)$ be a Riemannian manifold, $K\Subset M$ a bounded Lipschitz domain and $1<p<\infty$.
There exists $u\in A_p(K)$ with $\mathcal R_K(u)=(p-1)\Lambda_{p,M}(K)^p$. Consequently, by
Proposition \ref{prop:critical-attainment}, the functional $J_{\Lambda_{p,M}(K),p}$ has a nonzero global minimizer.
\end{proposition}

\begin{proof}
Fix $s>p$, with $s=Np/(N-p)$ if $p<N$ and $s=2p$ if $p\ge N$. Since $K$ is relatively compact,
the Sobolev and Poincaré inequalities in finitely many charts give a constant $C_S$ such that
$\|w\|_{L^s(K)}\le C_S\|\nabla_gw\|_{L^p(K)}$ for all $w\in W^{1,p}_0(K)$.

\emph{Lower bound on the plateaus.} For $u\in A_p(K)$, put $v=\Vol_g(\{u=1\})$. Since
$0\le u\le1$,
\[
v\le\int_Ku^s\dd V_g\le C_S^s\Big(\int_K|\nabla_gu|_g^p\dd V_g\Big)^{s/p},
\qquad\text{hence}\qquad
\mathcal R_K(u)\ge C_S^{-p}\,v^{\,p/s-1}.
\]
Because $p/s<1$, the right-hand side tends to $+\infty$ as $v\to0$. Let $(u_j)\subset A_p(K)$ be
a minimizing sequence and $v_j=\Vol_g(\{u_j=1\})$. Since $\mathcal R_K(u_j)$ is bounded, there is
$v_0>0$ with $v_j\ge v_0$ for all $j$.

\emph{Compactness.} Then $\int_K|\nabla u_j|^p=\mathcal R_K(u_j)\,v_j\le\mathcal R_K(u_j)\Vol_g(K)$
is bounded, so $(u_j)$ is bounded in $W^{1,p}_0(K)$. Up to a subsequence, $u_j\rightharpoonup u$
weakly in $W^{1,p}_0(K)$, $u_j\to u$ in $L^1(K)$ by the Rellich--Kondrachov theorem, and
$v_j\to\bar v\ge v_0$. Clearly $0\le u\le1$. For $t\in(0,1)$, the function $u_j-u$ is at least
$1-t$ on $\{u_j=1\}\setminus\{u>t\}$, so
\[
\Vol_g\big(\{u_j=1\}\setminus\{u>t\}\big)\le\frac{\|u_j-u\|_{L^1(K)}}{1-t}\longrightarrow0 .
\]
Hence $\Vol_g(\{u>t\})\ge\bar v$ for every $t\in(0,1)$. Since
$\{u=1\}=\bigcap_{t<1}\{u>t\}$, letting $t\uparrow1$ gives $\Vol_g(\{u=1\})\ge\bar v>0$, so
$u\in A_p(K)$. By weak lower semicontinuity,
\[
\mathcal R_K(u)\le\frac{\liminf_j\int_K|\nabla u_j|^p}{\bar v}=\lim_j\mathcal R_K(u_j)\,\frac{v_j}{\bar v}
=(p-1)\Lambda_{p,M}(K)^p ,
\]
and the reverse inequality holds by Proposition \ref{prop:threshold}. The last assertion follows from the
equivalence (1)$\Leftrightarrow$(2) in Proposition \ref{prop:critical-attainment}.
\end{proof}

\begin{remark}
Proposition~\ref{prop:attain} shows that at the critical value $\lambda=\Lambda_{p,M}(K)$ the
functional $J_{\lambda,p}$ always has a global minimizer other than $u=0$. It would be
interesting to know whether this minimizer is unique, and whether $v=1-u$ is a weak solution of
the Bernoulli problem in the sense recalled in  the Introduction.
%Section~\ref{sec:pos}.
\end{remark}

\subsection{Two lemmas}

\begin{lemma}[Strict monotonicity in the radius]\label{lem:mono}
Let $\kappa\in\R$ and $0<s<R$, with $R<\pi/\sqrt\kappa$ if $\kappa>0$. Then
$\Lambda_{p,\kappa}(B^\kappa_s)>\Lambda_{p,\kappa}(B^\kappa_R)$.
\end{lemma}

\begin{proof}
By Theorem~4.3 there is $\rho\in(0,s)$ with
$(p-1)\Lambda_{p,\kappa}(B^\kappa_s)^p=\mathrm{Cap}_{p,\kappa}(B^\kappa_\rho,B^\kappa_s)/V_\kappa(\rho)$.
By \eqref{eq:4-4}, $\mathrm{Cap}_{p,\kappa}(B^\kappa_\rho,B^\kappa_r)
=\omega_{N-1}\big(\int_\rho^rS_\kappa(t)^{-q}\dd t\big)^{1-p}$, with $q=(N-1)/(p-1)$, is strictly
decreasing in $r$ since $1-p<0$. Hence
\[
(p-1)\Lambda_{p,\kappa}(B^\kappa_R)^p\le\frac{\mathrm{Cap}_{p,\kappa}(B^\kappa_\rho,B^\kappa_R)}{V_\kappa(\rho)}
<\frac{\mathrm{Cap}_{p,\kappa}(B^\kappa_\rho,B^\kappa_s)}{V_\kappa(\rho)}=(p-1)\Lambda_{p,\kappa}(B^\kappa_s)^p .\qedhere
\]
\end{proof}

\begin{lemma}[Equality in the rearrangement inequality]\label{lem:eqps}
In the setting of Lemma~3.2 (respectively of Corollary~6.1, with normalized measures and
profiles), let $u\in W^{1,p}_0(K)$ with $0\le u\le1$ and $\Vol_g(\{u=1\})>0$, and assume that
\[
\int_K|\nabla_gu|_g^p\dd V_g=\int_{B}|\nabla u^\star|^p\dd V_{g_\kappa},
\]
where $B$ is the model ball and $u^\star$ the rearrangement. Then, for almost every $t\in(0,1)$,
the set $E_t=\{u>t\}$ satisfies $\Per_g(E_t)=I_\kappa\big(\Vol_g(E_t)\big)$. In particular $E_t$
is an isoperimetric region of $M$ and $I_M(\Vol_g(E_t))=I_\kappa(\Vol_g(E_t))$.
\end{lemma}

\begin{proof}
Let $\mu(t)=\Vol_g(\{u>t\})$; this function is nonincreasing, hence differentiable with
$\mu'(t)\in(-\infty,0]$ for a.e.\ $t$. Let $u$ also denote its precise representative, and put
$\Gamma_t=u^{-1}(t)\cap\{\nabla u\neq0\}$. By the coarea formula for Sobolev functions
\cite[3.2.12]{Fed69}, \cite[Section~2]{BZ88}, used in charts, the following hold for a.e.\
$t\in(0,1)$:
\begin{gather*}
\Per_g(E_t)\le\mathcal H^{N-1}(\Gamma_t),\qquad
-\mu'(t)\ \ge\ \int_{\Gamma_t}|\nabla u|^{-1}\dd\mathcal H^{N-1},\\
\int_K|\nabla u|^p=\int_0^1\!\!\int_{\Gamma_t}|\nabla u|^{p-1}\dd\mathcal H^{N-1}\dd t .
\end{gather*}
(The first relation is in fact an equality.) Hölder's inequality on $\Gamma_t$ gives
$\mathcal H^{N-1}(\Gamma_t)^p\le\int_{\Gamma_t}|\nabla u|^{p-1}\big(\int_{\Gamma_t}|\nabla u|^{-1}\big)^{p-1}$.
%Since $1-p<0$, and since $\Per_g(E)\ge I_M(\Vol_g(E))\ge I_\kappa(\Vol_g(E))$ by the definition
(%3.1) and 
%by (3.3) 
%(respectively by (6.2))
Since $1-p<0$, and since
$\Per_g(E)\ge I_M(\Vol_g(E))\ge I_\kappa(\Vol_g(E))$ by the definition \eqref{eq:3-1} and by
\eqref{eq:3-3} (respectively by \eqref{eq:6-2}), we obtain 
\begin{equation}\label{eq:chain}
\int_K|\nabla u|^p\ \ge\ \int_0^1\Per_g(E_t)^p\,(-\mu'(t))^{1-p}\dd t
\ \ge\ \int_0^1I_\kappa(\mu(t))^p\,(-\mu'(t))^{1-p}\dd t ,
\end{equation}
with the convention $0^{1-p}=+\infty$.

We claim that the last integral equals $\int_B|\nabla u^\star|^p$. Write
$u^\star=\varphi(d_{g_\kappa}(o,\cdot))$, with $\varphi$ nonincreasing and locally absolutely
continuous on $(0,R)$ (Lemma~3.1). For every $t\in(0,1)$ the set $\{u^\star>t\}$ is a concentric
ball of volume $\mu(t)$, whose boundary sphere has area $I_\kappa(\mu(t))$ by the model-ball realization
property,
 %by the model-ball realization hypothesis
  and $|\nabla u^\star|$ is constant on $\{u^\star=t\}$. Let
$A=\{r:\varphi'(r)=0,\ 0<\varphi(r)<1\}$. Since $\varphi$ is absolutely continuous,
$|\varphi(A)|\le\int_A|\varphi'|=0$. By the coarea formula, $\mu$ is the sum of the absolutely
continuous function $t\mapsto\int_t^1\int_{\{u^\star=\tau\}}|\nabla u^\star|^{-1}\dd\mathcal H^{N-1}\dd\tau$ and
of the distribution function of a measure carried by $\varphi(A)$, together with a constant on
$(0,1)$. The second part has zero derivative almost everywhere. Hence
$-\mu'(t)=\int_{\{u^\star=t\}}|\nabla u^\star|^{-1}$ for a.e.\ $t$; Hölder's inequality is an
equality on spheres, and the coarea formula for $u^\star$ proves the claim.

The two ends of \eqref{eq:chain} are therefore equal, which yields
\[
\int_0^1\Big(\Per_g(E_t)^p-I_\kappa(\mu(t))^p\Big)(-\mu'(t))^{1-p}\dd t=0 ,
\]
with a nonnegative integrand. For $t\in(0,1)$ we have
$\mu(t)\in[\Vol_g(\{u=1\}),V]$, so $I_\kappa(\mu(t))>0$; since $\int_B|\nabla u^\star|^p<\infty$,
this forces $-\mu'(t)>0$ for a.e.\ $t$. Hence $(-\mu'(t))^{1-p}\in(0,\infty)$ a.e., and
$\Per_g(E_t)=I_\kappa(\mu(t))$ for a.e.\ $t\in(0,1)$. The last assertion follows from
$I_\kappa(\mu(t))\le I_M(\mu(t))\le\Per_g(E_t)$.
\end{proof}

\subsection{Rigidity statements}

\begin{theorem}[Equality in Theorem~3.4]\label{thm:eq34}
Under the hypotheses of Theorem~3.4, assume that
$\Lambda_{p,M}(K)=\Lambda_{p,\kappa}(B_R^\kappa)$. Then
$\Per_g(K)=I_\kappa(\Vol_g(K))=I_M(\Vol_g(K))$: the domain $K$ is an isoperimetric region of $M$
whose perimeter equals that of the model ball. In particular, if $M=M_\kappa^N$ is the simply
connected space form of curvature $\kappa$, then equality holds in \eqref{eq:3-10} if and only
if $K$ is a geodesic ball of radius $R$. For $\kappa=0$ this recovers, for the variational
threshold, the rigidity part of the theorem of Daners and Kawohl \cite{DanersKawohl2010}.
\end{theorem}

\begin{proof}
\emph{Step 1: an optimal function.} Let $u$ be given by Proposition~\ref{prop:attain}. By
Lemma~3.2 and \eqref{eq:3-13}, $u^\star\in A_p(B^\kappa_R)$ and $\mathcal R_{B}(u^\star)\le\mathcal R_K(u)$,
where $B=B^\kappa_R$. Hence
\[
(p-1)\Lambda_{p,\kappa}(B^\kappa_R)^p\le\mathcal R_B(u^\star)\le\mathcal R_K(u)
=(p-1)\Lambda_{p,M}(K)^p=(p-1)\Lambda_{p,\kappa}(B^\kappa_R)^p ,
\]
so all these quantities are equal. Since the plateaus of $u$ and $u^\star$ have the same volume,
$\int_K|\nabla u|^p=\int_B|\nabla u^\star|^p$.

\emph{Step 2: $u$ is positive on $K$.} Suppose that $\Vol_g(\{u>0\})<V$. Then
$\{u^\star>0\}=B^\kappa_s$ for some $s<R$. Since $u^\star\in W^{1,p}_0(B)$ vanishes outside the
smooth ball $B^\kappa_s$, we have $u^\star\in W^{1,p}_0(B^\kappa_s)$, hence
$u^\star\in A_p(B^\kappa_s)$ and, by Lemma~\ref{lem:mono},
$\mathcal R_B(u^\star)\ge(p-1)\Lambda_{p,\kappa}(B^\kappa_s)^p>(p-1)\Lambda_{p,\kappa}(B^\kappa_R)^p$,
a contradiction. Hence $\Vol_g(\{u>0\})=V$; since $u=0$ a.e.\ outside $K$, $\{u>0\}=K$ up to a
null set.

\emph{Step 3: isoperimetric level sets.} By Lemma~\ref{lem:eqps},
$\Per_g(E_t)=I_\kappa(\Vol_g(E_t))$ for a.e.\ $t\in(0,1)$.

\emph{Step 4: passage to $K$.} Along a sequence $t_k\downarrow0$ of such levels,
$\mathbf 1_{E_{t_k}}\to\mathbf 1_{\{u>0\}}=\mathbf 1_K$ in $L^1$. By lower semicontinuity of the
perimeter and continuity of $I_\kappa$,
\[
\Per_g(K)\le\liminf_k\Per_g(E_{t_k})=\lim_kI_\kappa(\Vol_g(E_{t_k}))=I_\kappa(V)\le I_M(V)\le\Per_g(K),
\]
where the last two inequalities are \eqref{eq:3-3} and \eqref{eq:3-1}. This proves the first claim.

\emph{Step 5: space forms.} If $K$ is a geodesic ball of radius
$R$ in $M_\kappa^N$, equality holds by Section~8.1. Conversely, if $M=M^N_\kappa$, the isoperimetric regions are the geodesic balls,
up to null sets \cite[Chapter~VI]{Chavel}, \cite{BZ88b}; hence $|K\triangle B'|=0$ for a geodesic
ball $B'$ of radius $R$. Since $K$ is open, $K\setminus\overline{B'}$ is an open null set, hence
empty, so $K\subset\operatorname{int}\overline{B'}=B'$. If $x\in B'\setminus K$, then $x\in\partial K$:
otherwise a neighbourhood of $x$ in $B'$ would miss $K$, which has full measure in $B'$. Since
$K$ is a Lipschitz domain, $M\setminus K$ contains a truncated open cone with vertex $x$, whose
intersection with $B'$ has positive volume. This contradicts $|B'\setminus K|=0$. Hence $K=B'$.
\end{proof}

\begin{theorem}[Rigidity in Corollary~6.1]\label{thm:rigLG}
Let $(M^N,g)$ be closed and connected with $\mathrm{Ric}_g\ge(N-1)\kappa g$, $\kappa>0$, and let
$K\subset M$ be a Lipschitz domain with $\theta=\Vol_g(K)/\Vol_g(M)\in(0,1)$. Then equality holds
in \eqref{eq:6-3},
\[
\Lambda_{p,M}(K)=\Lambda_{p,\kappa}(B^\kappa_{r_\theta}),
\]
if and only if $M$ is isometric to the round sphere $S^N_\kappa$ and $K$ is a geodesic ball of
radius $r_\theta$.
\end{theorem}

\begin{proof}
If $M=S^N_\kappa$ and $K$ is a geodesic ball, equality holds by Example~6.3. Conversely, assume
equality. Steps 1--3 of the proof of Theorem~\ref{thm:eq34}, carried out with the normalized
measures $\mu_M,\mu_\kappa$, the normalized profiles and the rearrangement of Corollary~6.1, give
a function $u$ with $\{u>0\}=K$ up to a null set and, for a.e.\ $t\in(0,1)$,
\[
\frac{\Per_g(E_t)}{\Vol_g(M)}=I_\kappa\big(\mu_M(E_t)\big)=I_M\big(\mu_M(E_t)\big),
\qquad \mu_M(E_t)\in\big[\mu_M(\{u=1\}),\theta\big]\subset(0,1).
\]
Thus the normalized profile of $M$ coincides with that of $S^N_\kappa$ at some volume fraction
$v\in(0,1)$, and $E_t$ is an isoperimetric region realizing it. Equality in the Lévy--Gromov
inequality at one volume fraction forces $M$ to be isometric to $S^N_\kappa$
\cite[Appendix~C]{Gro99}; see \cite{CavallettiMondino} for the same rigidity on metric measure spaces. Step~4
then gives $\Per_g(K)/\Vol_g(M)=I_\kappa(\theta)$, so $K$ is an isoperimetric region of the round
sphere, and Step~5 shows that $K$ is a geodesic ball of radius $r_\theta$.
\end{proof}

\begin{remark}[Why not quote the rigidity of the Pólya--Szegő inequality directly]
Mondino and Semola proved a rigidity theorem for the Pólya--Szegő inequality on
$\mathrm{RCD}(N-1,N)$ spaces \cite[Theorem~5.1]{MS20}, under a nondegeneracy assumption on the
gradient ($|\nabla u|\neq0$ almost everywhere on the support of $u$). This assumption fails for
the functions relevant here, which have a plateau $\{u=1\}$ of positive volume. The proof above
avoids it by using only the levels $t\in(0,1)$ and the rigidity of the isoperimetric
inequality at a single volume.
\end{remark}

%\begin{remark}[Equality in Corollary~\ref{cor:avr}]\label{rem:avrrig}
%The same proof, with Lemma~\ref{lem:ps} in place of Lemma~3.2 of the article and Euclidean model balls, shows
%that equality in Corollary~\ref{cor:avr} forces
%\[
%\Per_g(K)=\mbox{AVR}_g^{1/N}\,I_0(\Vol_g(K)),
%\]
%that is, equality in the sharp isoperimetric inequality of Lemma~\ref{lem:iso}. By the rigidity theorem of Cavalletti and Manini \cite{CMa24}, equality in
%that inequality can be attained only by metric balls; on a smooth manifold this is expected to
%force $\mbox{AVR}_g=1$, hence $M$ isometric to $\R^N$ and $K$ a ball. The precise statement to be
%invoked for Lipschitz domains in smooth manifolds should be checked against \cite{CMa24} before
%this is stated as a theorem. In the non-smooth setting, the metric cones of the previous remark
%are equality cases, consistently with the computation given there.
%\end{remark}
\begin{theorem}[Rigidity in Corollary~6.8]
Let $(M,g)$, $K$ and $B$ be as in Corollary~6.8, and let $\theta=\mbox{AVR}_g$. Then
\[
\Lambda_{p,M}(K)=\theta^{1/N}\,\Lambda_{p,0}(B)
\]
if and only if $M$ is isometric to $\R^N$ and $K$ is a ball.
\end{theorem}

\begin{proof}
If $M=\R^N$, then $\theta=1$, and if $K$ is a ball, $K$ is isometric to $B$ and equality holds.
Conversely, assume equality. For $u\in A_p(K)$ we use the notation $\mathcal R_K(u)$ of Section~7.

\emph{Step 1: an optimal function.} Let $u$ be given by Proposition~7.1, so that
$\mathcal R_K(u)=(p-1)\Lambda_{p,M}(K)^p$. By Lemma~6.7, its Schwarz rearrangement $u^\sharp$
belongs to $W^{1,p}_0(B)$, and, as in \eqref{eq:3-13}, $|\{u^\sharp=1\}|=\Vol_g(\{u=1\})>0$, so
that $u^\sharp\in A_p(B)$ and $\theta^{p/N}\mathcal R_B(u^\sharp)\le\mathcal R_K(u)$. Hence
\[
\theta^{p/N}(p-1)\Lambda_{p,0}(B)^p\le\theta^{p/N}\mathcal R_B(u^\sharp)\le\mathcal R_K(u)
=(p-1)\Lambda_{p,M}(K)^p=\theta^{p/N}(p-1)\Lambda_{p,0}(B)^p ,
\]
so that all these quantities are equal. In particular
$\mathcal R_B(u^\sharp)=(p-1)\Lambda_{p,0}(B)^p$ and
$\int_K|\nabla_gu|_g^p\dd V_g=\theta^{p/N}\int_B|\nabla u^\sharp|^p\dd x$.

\emph{Step 2: $u$ is positive on $K$.} If $\Vol_g(\{u>0\})<\Vol_g(K)$, then
$\{u^\sharp>0\}$ is a ball $B_s$ of radius $s<R^e_K$, and $u^\sharp\in A_p(B_s)$. By Lemma~7.3
with $\kappa=0$,
$\mathcal R_B(u^\sharp)\ge(p-1)\Lambda_{p,0}(B_s)^p>(p-1)\Lambda_{p,0}(B)^p$, a contradiction.
Hence $\{u>0\}=K$ up to a null set.

\emph{Step 3: isoperimetric level sets.} Let $E_t=\{u>t\}$ and $\mu(t)=\Vol_g(E_t)$. The coarea
and Hölder argument of Lemma~7.4, with Lemma~6.6 in place of \eqref{eq:3-3}, gives
\begin{align*}
\int_K|\nabla_gu|_g^p\dd V_g&\ge\int_0^1\Per_g(E_t)^p(-\mu'(t))^{1-p}\dd t\\
&\ge\theta^{p/N}\int_0^1I_0(\mu(t))^p(-\mu'(t))^{1-p}\dd t
=\theta^{p/N}\int_B|\nabla u^\sharp|^p\dd x ,
\end{align*}
where the last identity is proved exactly as in Lemma~7.4 with $\kappa=0$. By Step~1 the two
ends are equal. As in Lemma~7.4, $I_0(\mu(t))>0$ forces $-\mu'(t)>0$ for a.e.\ $t\in(0,1)$, and
we obtain $\Per_g(E_t)=\theta^{1/N}I_0(\mu(t))$ for a.e.\ $t\in(0,1)$.

\emph{Step 4: equality for $K$.} Along a sequence $t_k\downarrow0$ of such levels,
$\mathbf 1_{E_{t_k}}\to\mathbf 1_K$ in $L^1$. By lower semicontinuity of the perimeter and
Lemma~6.6,
\[
\Per_g(K)\le\liminf_k\Per_g(E_{t_k})=\theta^{1/N}I_0(\Vol_g(K))\le\Per_g(K).
\]
Thus $K$ realizes equality in the sharp isoperimetric inequality of Lemma~6.6.

\emph{Step 5: rigidity.} The metric measure space $(M,d_g,\Vol_g)$ is a non-collapsed
$\mathrm{RCD}(0,N)$ space with Euclidean volume growth. By the rigidity part of the sharp
isoperimetric inequality \cite[Theorem~1.1]{APPS23}, $M$ is isometric to a Euclidean metric
measure cone with tip $o$, and $K$ coincides, up to a null set, with a ball $B_r(o)$. In a
Euclidean metric measure cone, $\Vol_g(B_s(o))=\theta\beta_Ns^N$ for every $s>0$. Since $M$ is
smooth at $o$, $\Vol_g(B_s(o))/(\beta_Ns^N)\to1$ as $s\to0$. Hence $\theta=1$, and $M$ is
isometric to $\R^N$ by the equality case of the Bishop--Gromov inequality. Finally, $K$ is an
open Lipschitz set which coincides with the ball $B_r(o)$ up to a null set; as in Step~5 of the
proof of Theorem~7.5, $K=B_r(o)$.

When $\partial K$ is smooth, Step~5 can be replaced by the equality case of Brendle's sharp
Sobolev inequality \cite[Theorem~1.2]{Bre21}, applied with $f\equiv1$ and $D=\overline K$.
\end{proof}

\begin{remark}[Non-smooth spaces]
Steps~1--4 use only the attainment of the threshold, the Pólya--Szegő inequality with constant
and the coarea formula. If Corollary~6.8 extends to non-collapsed $\mathrm{RCD}(0,N)$ spaces with
Euclidean volume growth, the same argument and \cite[Theorem~1.1]{APPS23} show that equality
holds exactly when the space is a Euclidean metric measure cone and $K$ a ball centred at one of
its tips. This is consistent with Remark~6.10.
\end{remark}

%===========================================================
%===========================================================
\section{Examples}
\label{sec:conclusion}
%===========================================================

\subsection{Equality in the model case}

If $M=M_\kappa^N$ and $K=B_R^\kappa$, Theorem~\ref{thm:main} is an equality:
\[
\Lambda_{p,M_\kappa^N}(B_R^\kappa)
=
\Lambda_{p,\kappa}(B_R^\kappa).
\]
The value is given by the explicit formula \eqref{eq:4-7}. In the Euclidean case, it reduces to \eqref{eq:4-13}, and for $N=3,p=2$ it is $9/(2R)$.

\subsection{Further problems}

\noindent\textit{Notation.} Here and throughout, \(\operatorname{Cap}_{p,\kappa}(E,B_R^\kappa)\) denotes the relative $p$-capacity of $E$ with respect to $B_R^\kappa$, as defined in Section~2.

The results proved above lead to the following two precise problems.

\begin{enumerate}
\item[\textup{(P1)}] \textbf{Small-radius expansion of the variational threshold.}
Let $(M^N,g)$ be a smooth Riemannian manifold, let $x_0\in M$, and
let $1<p<\infty$. Determine the asymptotic behavior, as $R\downarrow0$,
of
\[
\Lambda_{p,M}(B_R(x_0)).
\]
In particular, decide whether there exists a constant
$\gamma_{N,p}$, depending only on $N$ and $p$, such that
\[
\Lambda_{p,M}(B_R(x_0))
=
\frac{C_{N,p}}{R}
+
\gamma_{N,p}\,\operatorname{Scal}_g(x_0)\,R
+o(R),
\qquad R\downarrow0,
\]
where
\[
C_{N,p}=\Lambda_{p,0}(B_1)
=\frac{N}{(p-1)x_{N,p}}
\]
is the Euclidean unit-ball threshold given by \eqref{eq:4-13}.
If this expansion does not hold in this form, identify the curvature
invariants and the correct order of the first geometric correction.
A proof must control the optimizing positivity sets in the capacity
characterization uniformly as $R\downarrow0$; an expansion of the
capacity of a prescribed geodesic ball alone does not establish the
claimed expansion for the variational threshold.

\medskip
%\noindent\textit{Notation.} Here and throughout, \(\operatorname{Cap}_{p,\kappa}(E,B_R^\kappa)\) denotes the relative $p$-capacity of $E$ with respect to $B_R^\kappa$, as defined in Section~2.
\item[\textup{(P2)}] \textbf{Quantitative stability of the model capacity minimizer.}
Fix $N\geq2$, $1<p<\infty$, a space-form parameter $\kappa$, and a
radius $R$ for which the model ball $B_R^\kappa$ is defined. Set
\[
q=\frac{N-1}{p-1},\qquad
F_\kappa(\rho)=
\frac{\operatorname{Cap}_{p,\kappa}
(B_\rho^\kappa,B_R^\kappa)}
{(p-1)\operatorname{Vol}_\kappa(B_\rho^\kappa)},
\quad 0<\rho<R.
\]
The explicit formula in Section~4 determines the minimum of
$F_\kappa$. Determine whether its minimizing radius $\rho_*$ is
unique and nondegenerate, and, whenever it is, whether there exists a
constant $c=c(N,p,\kappa,R)>0$ such that every $p$-quasi-closed set
$E\subset B_R^\kappa$ with positive volume satisfies
\[
\frac{\operatorname{Cap}_{p,\kappa}
(E,B_R^\kappa)}
{(p-1)\operatorname{Vol}_\kappa(E)}
-\Lambda_{p,\kappa}(B_R^\kappa)^p
\geq
c\left(
\frac{\operatorname{Vol}_\kappa
(E\mathbin{\triangle}B_{\rho_*}^\kappa)}
{\operatorname{Vol}_\kappa(B_R^\kappa)}
\right)^2.
\]
Here $B_{\rho_*}^\kappa$ is the concentric minimizing ball and
$\triangle$ denotes symmetric difference. The uniqueness and
nondegeneracy of $\rho_*$, as well as the uniform estimate over all
admissible quasi-closed sets, are part of the problem and are not
assumed by the present paper. A proof would require a quantitative
rearrangement or isoperimetric inequality together with a coercive
estimate for the one-dimensional minimization problem.

 Let us mention that for
$\kappa=0$ the optimal radius is unique by \eqref{eq:4-12}.

\end{enumerate}

%\textbf{Acknowledgements and declarations}\\
%The author thanks Claude (Anthropic), an artificial intelligence assistant, for its help in checking
%the computations, implementing the numerical method and preparing the manuscript.
%The author takes full responsibility for the entire content of this article: statements, proofs,
%numerical results and conclusions.

\noindent\textbf{Declaration of generative AI and AI-assisted technologies in the manuscript
preparation process.} During the preparation of this work, the authors used Claude (Anthropic;
model Claude Opus~5.5) in order to check computations and proofs, to suggest and draft parts of
the text and of some proofs (in particular in Sections~2.1, 6.3 and~7), and to improve the language. After using this tool, the
authors reviewed, verified and edited all the content as needed and take full responsibility for
the content of the publication.

\noindent\textbf{Acknowledgements.} The authors thank Claude (Anthropic), an artificial
intelligence assistant, for its help in checking computations and drafting parts of the manuscript. The authors take full responsibility for
the entire content of this article: statements, proofs and conclusions.

\noindent\textbf{Funding.} This research received no specific grant from any funding agency. %/
%indicate the grant.]

\noindent\textbf{Competing interests.} The authors declare that they have no competing
interests.

\noindent\textbf{Data availability.} Data sharing is not applicable to this article, as no datasets were generated or analysed during the current study.

\end{document}